\documentclass[review,hidelinks,onefignum,onetabnum]{siamart250211}

\usepackage{lipsum}
\usepackage{amsfonts}
\usepackage{graphicx}
\usepackage{epstopdf}
\usepackage{algorithmic}

\usepackage{subfig}

\usepackage[skip=4pt]{caption}        % caption 与图之间
\enlargethispage{\baselineskip}

\usepackage{etoolbox}  % 方便一次性修改
\ifpdf
  \DeclareGraphicsExtensions{.eps,.pdf,.png,.jpg}
\else
  \DeclareGraphicsExtensions{.eps}
\fi

\newsiamremark{remark}{Remark}
\newsiamremark{hypothesis}{Hypothesis}
\crefname{hypothesis}{Hypothesis}{Hypotheses}
\newsiamthm{claim}{Claim}
\newsiamremark{fact}{Fact}
\crefname{fact}{Fact}{Facts}

\headers{A domain decomposition approach to pore-network modeling}{S. Sun, Z. Wang, L. Zhang and J. Zhao}

\title{A domain decomposition approach to pore-network modeling of porous media flow\thanks{Submitted to the editors DATE.
}}

\author{Shuyu Sun\thanks{School of Mathematical Sciences, Key Laboratory of Intelligent Computing and Applications (Ministry of Education), Tongji University, Shanghai 200092, China  
  (\email{suns@tongji.edu.cn},
  \email{1954013@tongji.edu.cn}, \email{22210@tongji.edu.cn},
\email{jijingzhao@tongji.edu.cn}
  ).}
\and Zhangchengrui Wang\footnotemark[2]  
\and Lei Zhang\footnotemark[2]
\thanks{Corresponding author, 
\email{22210@tongji.edu.cn}.}
\and Jijing Zhao\footnotemark[2]
}

\usepackage{amsopn}



\begin{document}

\maketitle

\nolinenumbers

% REQUIRED
\begin{abstract}
We propose a domain-decomposition pore-network method (DD-PNM) for modeling single-phase Stokes flow in porous media. The method combines the accuracy of finite-element discretizations on body-fitted meshes within pore subdomains with a sparse global coupling enforced through interface unknowns. 
Local Dirichlet-to-Neumann operators are precomputed from finite-element solutions for each pore subdomain, enabling a global Schur complement system defined solely on internal interfaces.
Rigorous mathematical analysis establishes solvability and discrete mass conservation of the global system. 
Moreover, we constructively recover classical pore-network models by fitting half-throat conductivities to local Dirichlet-to-Neumann maps, providing a principled bridge between mesh-based and network-based frameworks.
Numerical results are presented to demonstrate the validity and effectiveness of the overall methodology.
\end{abstract}

% REQUIRED

% REQUIRED

\section{Introduction}
\label{sec:intro}

Accurate pore-scale simulations are essential to predict and design transport phenomena in porous materials across energy and environmental systems. Representative applications include hydrocarbon recovery and enhanced oil recovery, subsurface hydrogen and CO$_2$ storage and leakage risk assessment, digital rock analysis for reservoir characterization, and water management in fuel-cell gas diffusion layers and battery electrodes \cite{blunt2013pore,cui2022pore,gostick2016openpnm,cooper2017simulated,shan2022super,tranter2016pore}. 
These applications demand models that (i) respect pore‑scale physics and geometry, (ii) expose macroscopic observables with quantified errors, and (iii) remain computationally viable on large image volumes and across parameter sweeps.

Current approaches to pore-scale simulations fall broadly into two families. 
(i) \emph{Direct modeling} solves the Stokes or Navier-Stokes equations on voxelized micro-CT images using finite volume/element (FV/FE) solvers or Lattice Boltzmann (LB) methods \cite{bijeljic2013predictions,yang2019fully,zhang2019coupled}. Direct modeling can resolve complex geometries and provide full field information, but it entails very large degrees of freedom, so scaling to centimeter-scale samples or performing many-query studies (e.g., parameter estimation, uncertainty quantification) is challenging \cite{icardi2014pore,mehmani2025multiscale,mcclure2021lbpm}. 
(ii) \emph{Classical pore-network modeling} (C-PNM) \cite{blunt2013pore,xiong2016review,cui2022pore} replaces the voxel image by a graph of pores (nodes) connected by throats (edges), and solves the flow in the network by using one-dimensional (1-D) analytical solution in each edge together with mass conservation at each node. 
% and solves one-dimensional (1-D) flow in each edge with mass conservation at nodes.
% This affords orders-of-magnitude speed-ups, but at the cost of strong modeling assumptions. 
 This affords orders-of-magnitude speed-ups, but at the cost of strong assumptions of model reduction.

Since Fatt’s seminal work on PNM \cite{fatt1956network,blunt2013pore,xiong2016review,cui2022pore}, substantial progress has been made in efficient image-driven network extraction, in refining physics-based closures and parameterization (e.g., conductivity), and in broadening applicability to complex wettability, multiphase flow, and irregular geometries. Maximal ball methods identify pore bodies as local maxima of the distance map and cluster inscribed spheres to define throats \cite{silin2006pore,dong2009pore}. 
Medial axis approaches \cite{lindquist1996medial,lindquist1999investigating,prodanovic20073d,jiang2007efficient,yi2017pore}
reduce the void space to a centerline graph whose junctions and links define pores and throats. 
Morphology-based workflows use distance transforms, erosion-dilation, and watershed segmentation to delineate pores and throats \cite{rabbani2014automated,gostick2017versatile,
zhao2020simulation}. 
Comprehensive overviews of extraction strategies can be found in \cite{blunt2013pore,xiong2016review,cui2022pore}. More recently, advanced network extraction techniques move beyond pixel-based skeletons toward "pixel-free" or hybrid medial-axis approaches, thereby reducing discretization artifacts and truncation sensitivity \cite{liu2024new,zhang2024medial}. 

\begin{figure}[htbp!]
\centering
\includegraphics[width=7cm]{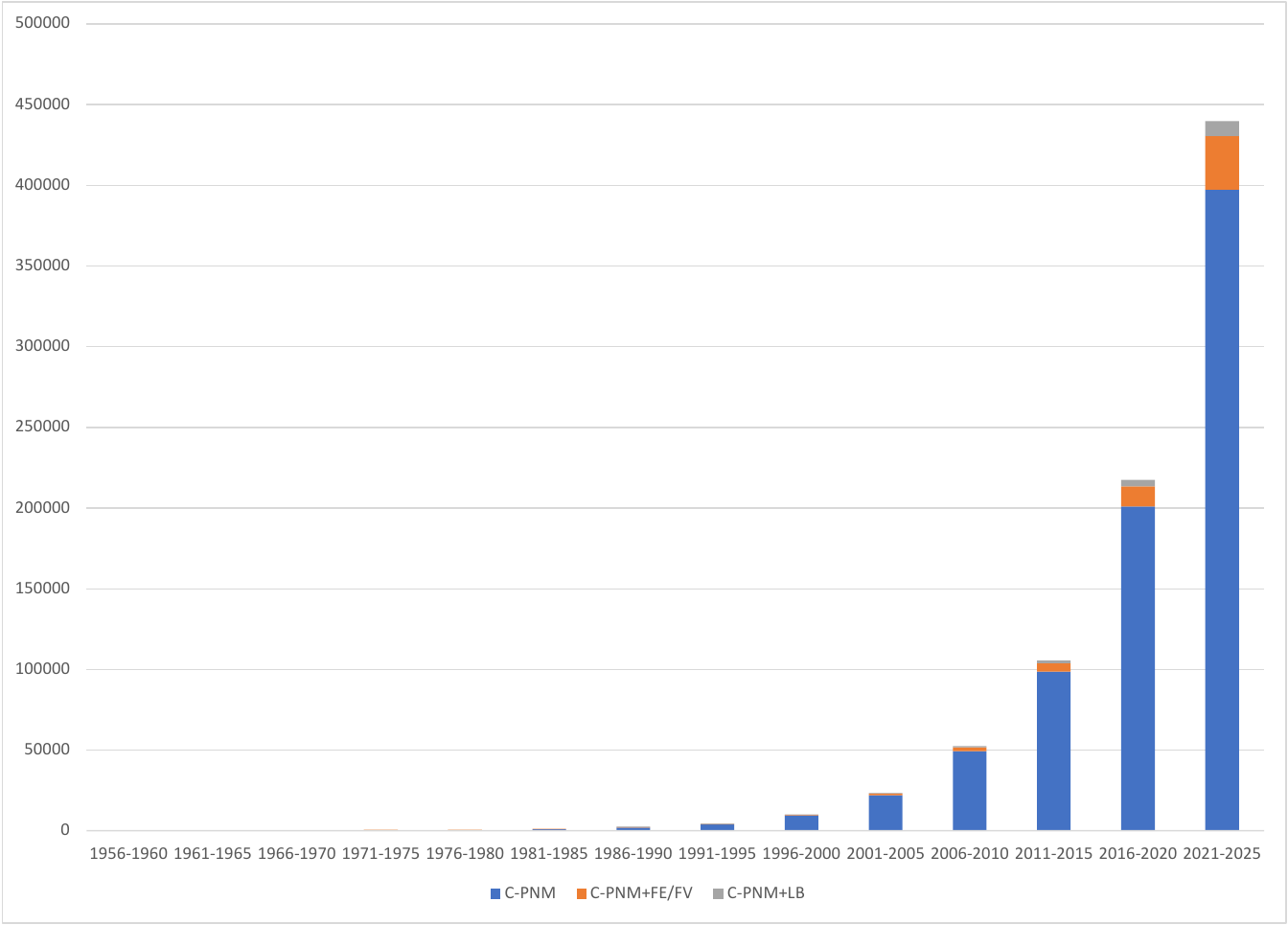}
\caption{Publication counts (every 5 years) for C-PNM, C-PNM+FE/FV, and C-PNM+LB (1956-2025); data source: Scopus \cite{scopus2024}.} 
\label{fig:pnm_trends}
\end{figure}

Beyond network extraction, model fidelity hinges on the physics assigned to the network and on how pore-scale parameters (e.g., conductivity, permeability)  are obtained. Classical one-dimensional Poiseuille-type closures with shape factors remain the baseline for specifying throat conductivities in C-PNM \cite{cui2022pore}. Generalized network models (GNM) enrich pore/throats with corners and shape descriptors to better  capture capillarity, wettability \cite{raeini2017generalized,raeini2018generalized}. 
%Figure \ref{fig:pnm_trends} highlights the rapid growth of PNM across diverse applications, together with a rising body of studies that benchmark or combine PNM with direct simulations using FE/FV/LBM. In these works, 
Increasingly, parameter values are obtained by regressing from local direct simulations (FV/FE/LB) on segmented subdomains, and by benchmarking network predictions against direct simulations \cite{prodanovic20073d,tansey2016pore,
raeini2018generalized,zhao2020simulation,mcclure2021lbpm,
kohanpur2022using,zhao2023pore}. 
Coupled formulations with C-PNM have also been explored: \cite{weishaupt2019efficient,weishaupt2020hybrid} formulate a fully monolithic model coupling Navier-Stokes flow in a free-flow region to a C-PNM in the porous region, and  \cite{fagbemi2020coupling} couples a multiphase C-PNM to an FE description of structural deformation for fluid-solid interaction in deformable media. 
Figure \ref{fig:pnm_trends} summarizes publication trends in pore-scale modeling (1956-2025), including C-PNM and studies pairing C-PNM with FE/FV/LB for calibration or region coupling. The figure highlights a notable growth in C-PNM studies over recent years and the increasing use of combined C-PNM and continuum/mesoscopic methods for enhanced modeling accuracy.

% have surged over the past decade, and studies that couple PNM with continuum or mesoscopic solvers (PNM+FEM, PNM+LBM) are increasingly common. 

%Most of these works implement region coupling or calibration between distinct solvers. In contrast, our approach provides a constructive bridge: we retain a mesh-based finite-element description inside pores, derive DtN-based interface operators, and assemble a sparse SPD interface system that can be reduced to C-PNM via fitted half-throat conductivities. 

%????The plot highlights the rapid growth of PNM and the increasing interest in couplings with continuum and mesoscopic solvers. Our study differs from such couplings by offering a bridge to mesh-based solvers: FE within pores, DtN-derived interface operators, and a sparse SPD interface system reducible to C-PNM. Data sources, keywords, and de-duplication rules are provided in ????. \\
%???? To situate our contribution, Figure \ref{fig:pnm_trends} charts publication trends in pore‑scale modeling (1956–2025), including PNM alone and studies combining PNM with FE/FV/LBM for calibration or region coupling. \\
%Figure \ref{fig:pnm_trends}. Annual publication counts for PNM, PNM+FEM, and PNM+LBM (1956–2025). The plot highlights rapid growth in PNM and increasing interest in PNM–continuum/mesoscopic combinations. Our study differs from such couplings by providing a constructive bridge: FE inside pores, DtN‑derived interface operators, and a sparse SPD interface system reducible to C‑PNM. Data sources, search queries, and de‑duplication rules are detailed in ????.

Despite recent advances and the wide use of C-PNM
for fast upscaling, several limitations persist: (i) geometrical idealization and 1-D closure laws: each throat is typically treated as a Poiseuille conduit with simple cross-section, which reduces fidelity in the presence of strong shear, corner/film flows, or highly irregular cross-sections \cite{blunt2013pore,liu2024new};
(ii) limited internal field resolution: C-PNM predicts fluxes and pressures at network elements but not detailed intra-pore fields required for, e.g., heat/mass transfer coupling or shear-dependent rheology \cite{xiong2016review,blunt2013pore};
(iii) sensitivity to segmentation and network extraction:  predictions depend on image processing and network extraction algorithms (maximal balls, medial axis, etc.), which can be sensitive to noise, grayscale thresholding, and resolution \cite{dong2009pore,lindquist1996medial,yi2017pore,
zhang2024medial,gostick2017versatile};
 (iv) limited error quantification and convergence theory: rigorous error estimates relative to mesh-based reference solutions, and a constructive bridge to FE/FV formulations remain limited \cite{xiong2016review,mehmani2018multiscale}. 
These limitations motivate a framework that preserves mesh-level fidelity within pores while retaining a sparse, network-like global coupling and providing an analyzable, constructive bridge to C-PNM. 

%Given these limitations of C-PNM, direct modeling is a natural alternative, but its computational cost remains prohibitive for large image volumes and parameter studies. 
%To bridge the gap between the high fidelity of direct modeling and its computational expense,

To this end, we adopt a domain-decomposition perspective that replaces a monolithic solve by local subproblems coupled through interface conditions, thereby reducing global cost while retaining accuracy  \cite{liu2002domain,rabbani2019pore,stavroulakis2009advances,quarteroni1999domain}. 
In this work, we tailor this paradigm to align with  pore‑network structure: 
non-overlapping subdomains are defined to coincide with pore bodies in C-PNM, and we enforce a single scalar normal-traction unknown per internal face, consistent with the network view of dominant normal exchange.

\textbf{Contributions.} We propose and analyze a domain decomposition pore-network method (DD-PNM) that combines mesh-based accuracy inside pores with a sparse network coupling across interfaces:
(i) formulation and analysis. Each pore body is an independent finite-element subdomain for Stokes flow, discretized using an unstructured mesh generated from the segmented void geometry;
for each pore subproblem, we precompute a local Dirichlet-to-Neumann (DtN) map that relates normal traction to normal flux; enforcing global mass conservation across all interfaces yields a sparse Schur complement system with size equal to the number of interfaces; 
 we establish solvability, symmetric positive-definiteness (SPD) of the Schur complement system, and both local and global discrete mass conservation. This formulation retains mesh-level fidelity within pores, is robust to irregular geometries, and naturally supports offline-online workflows (via precomputation of local operators) and parallel implementations; 
(ii) constructive bridge to C-PNM. We derive a reduction from DD-PNM to C-PNM by fitting half-throat conductivities directly from the local DtN responses, thus providing a mesh-informed, non-empirical pathway form pore-scale physics to network parameters \cite{raeini2017generalized,raeini2018generalized}; 
%This clarifies the relation between mesh-based physics and graph-based models. 
(iii) error estimate. We compare the proposed DD-PNM with a classical FE domain-decomposition formulation and  derive an error estimate of the DD-PNM.

 \textbf{Relation to previous works.} Our work is related to the pore-level multiscale method (PLMM) \cite{mehmani2018multiscale}, which constructs local response operators and eliminates fine-scale degrees of freedom, resulting in a sparse global coupling among interface unknowns. While PLMM uses finite-volume discretizations on Cartesian grids , our framework instead solve subdomain Stokes problems using body-fitted finite elements on unstructured meshes conforming to pore geometry. 
%Furthermore, we explicitly bridge the gap between mesh-based and graph-based pore-network models by deriving a constructive reduction from DD-PNM to C-PNM through fitting half-throat conductivities to local Dirichlet-to-Neumann (DtN) maps. 
In addition, the boundary conditions for the local problems in our framework differs from those used in PLMM. 
 We also provide a rigorous mathematical analysis, demonstrating solvability, symmetry, positive definiteness of the Schur complement, and both local and global discrete mass conservation of the assembled DD-PNM system. 
Numerous studies infer equivalent local parameters by performing direct simulations (FE/FV/LB) on segmented subdomains and benchmark C-PNM predictions against these high-fidelity solvers \cite{prodanovic20073d,tansey2016pore,raeini2018generalized,zhao2020simulation,mcclure2021lbpm,kohanpur2022using,zhao2023pore}. We extend beyond parameter calibration by building DtN operators for each pore subproblem and assembling a global interface system. Crucially, our approach admits a constructive reduction to C-PNM by fitting half-throat parameters directly to DtN responses, thereby providing a constructive bridge between mesh-based physics and graph-based models. 
Finally, monolithic or hybrid couplings have been proposed to couple  continuum/mesoscopic solvers with  pore-network solvers,
effectively leveraging the computational efficiency of network models  \cite{weishaupt2019efficient,weishaupt2020hybrid,fagbemi2020coupling}.  
These coupling approaches improve computational efficiency by explicitly dividing the domain into distinct regions, solving expensive fine-scale models (e.g., continuum or mesoscopic simulations) only where necessary, and using simplified pore-network models elsewhere.
 In contrast, our contribution integrates all pore-resolved finite-element subproblems via DtN-based interface operators into a unified, sparse global system. Efficiency is achieved as the global unknowns scale only with the number of internal interfaces, maintaining both computational scalability and mesh-level resolution inside pore geometries.

The outline of the paper is as follows. 
 Section 2 introduces the methodology, beginning with the C-PNM, followed by the reference finite-element formulation used for validation and local operator construction. We then derive the DD-PNM, detailing the  local Dirichlet to Neumann maps and global assembly. Additionally, we explain the constructive reduction from DD-PNM to C-PNM. Section 3 presents the mathematical analysis establishing solvability and discrete mass conservation properties. Numerical results validating the effectiveness of the DD-PNM are provided in Section 4. Finally, Section 5 concludes the paper and outlines directions for future research.

\section{Methodology}\label{sec:method}

\subsection{Model problem}
\label{sec:model_problem}

\begin{figure}[htbp!]
  \centering
%\subfloat[]{\includegraphics[width=0.36\textwidth]{figures/particle.eps}}
%~~
%\subfloat[]{\includegraphics[width=0.36\textwidth]{figures/PNM_model.eps}}
%\\
%\vspace*{-0.4cm}
%\subfloat[]{\includegraphics[width=0.36\textwidth]{figures/FEMregion.eps}}
%~~
%\subfloat[]{\includegraphics[width=0.36\textwidth]{figures/DDPNMmesh.eps}}
\includegraphics[width=1\textwidth]{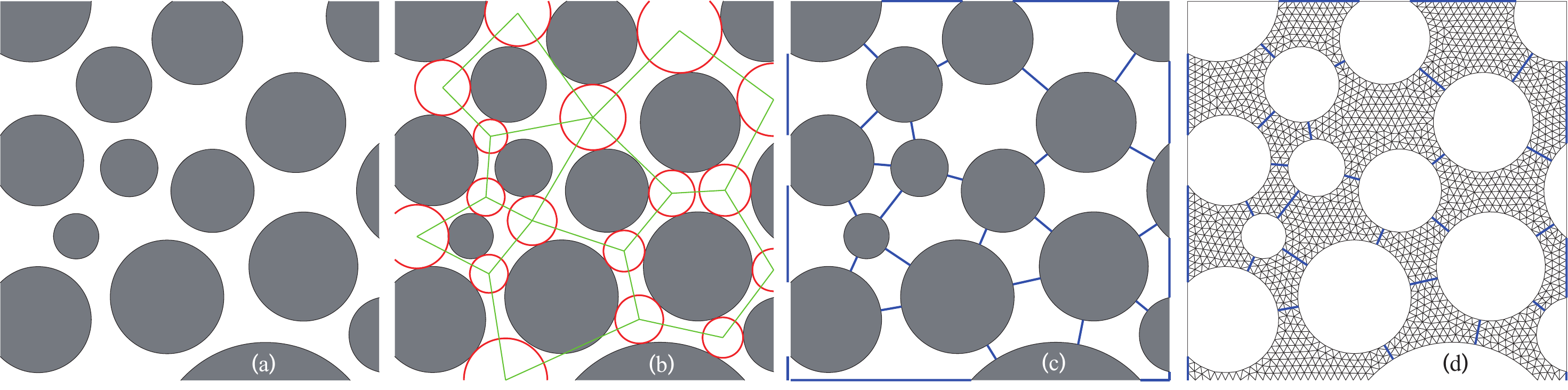}
  \caption{(a) computational geometry;  
           (b) maximal ball graph for C-PNM;  
           (c) subdomains $\{\Omega_i\}$ used by DD-PNM;  
           (d) body-fitted FE mesh.}
  \label{fig:model_pb}
\end{figure}

We consider steady, incompressible Stokes flow through a
two-dimensional porous medium.  
The synthetic micro-geometry, shown in Figure \ref{fig:model_pb}(a), is
obtained by packing circular solid particles (grey) inside a rectangular domain; the remaining void (white) represents the pore space.  
Throughout the present section, this simple yet
representative  layout serves as a running example for explaining each ingredient of the proposed methodology. 
The same geometry is revisited in Section \ref{sec:num_results}, where numerical results are produced and the performance of the method is quantitatively assessed.

Let $\Omega\subset\mathbb R^{2}$ denote the pore domain, and let
$\partial\Omega^{\text{wall}}$ be the solid-fluid interface.
A pressure drop $p_{\mathrm{in}}-p_{\mathrm{out}}>0$ is imposed
between the left ($\Gamma_{\mathrm{in}}$) and right
($\Gamma_{\mathrm{out}}$) boundaries, whereas the top ($\Gamma_{\rm top}$) and bottom ($\Gamma_{\rm bot}$) boundaries are treated as  no-slip
walls. We collect all homogeneous Dirichlet segments in the single set $\Gamma_{\mathrm{dir}} := \partial\Omega^{\text{wall}}\cup \Gamma_{\mathrm{top}} \cup \Gamma_{\mathrm{bot}}$. The governing equations read
\begin{equation}
\left\{
\begin{array}{lll}
-\nu \Delta \mathbf u + \nabla p = \mathbf 0, &&\text{in } \Omega, \\
\nabla\!\cdot\!\mathbf u = 0, & &
\text{in } \Omega, \\
\mathbf u = \mathbf 0,
&&\text{on } \Gamma_{\rm dir}, \\
\mathbf n\!\cdot\!\boldsymbol\sigma(\mathbf u,p)\mathbf n = -p_{\mathrm{in}},\; & \mathbf t\!\cdot\!\boldsymbol\sigma(\mathbf u,p)\mathbf n = 0, & \text{on }\Gamma_{\mathrm{in}},\\
\mathbf n\!\cdot\!\boldsymbol\sigma(\mathbf u,p)\mathbf n = -p_{\mathrm{out}},\; & \mathbf t\!\cdot\!\boldsymbol\sigma(\mathbf u,p)\mathbf n = 0, & \text{on }\Gamma_{\mathrm{out}}.\\
\end{array}
\right.
\label{eq:stokes-global}
\end{equation}
Here, $\mathbf u$ is the velocity field, $p$ is the pressure, 
$\nu$ is the kinematic viscosity, $\boldsymbol{\sigma}(\mathbf u,p):=-p\mathbf{I}+\nu\nabla\mathbf{u}$ is the Cauchy stress tensor, $\mathbf{n}$ is the outward unit normal, and $\mathbf t\perp \mathbf n$.

%The geometric segmentation required by each modelling strategy

%Maximal-ball skeletonisation partitions the pore space into
%$N$ non-overlapping, Lipschitz subdomains
%$\{\Omega_i\}_{i=1}^{N}$, each representing a pore body.  For every
%adjacent pair $(i,j)$ we denote by
%$\Gamma_{ij}:=\partial\Omega_i\cap\partial\Omega_j$ the median
%cross-section of their connecting throat.

\noindent
\textbf{Road map.}\;
With the geometry and governing equations fixed, we proceed in three steps: \\
1.~~\emph{Baseline:}  we briefly introduce the classical pore-network model
(C-PNM) 
(Section~\ref{sec:cpnm}).\\[2pt]
2.~~\emph{Reference:}  we present a full finite-element (FE) resolution 
serving both as a validation benchmark and as a building block for our proposed
method (Section~\ref{sec:fem}).\\[2pt]
3.~~\emph{Hybrid:}  we derive the domain-decomposition pore-network
method (DD-PNM), which combines FE accuracy with PNM scalability
(Sections~\ref{sec:ddpnm}).

\subsection{Classical pore–network modeling (C-PNM)}
\label{sec:cpnm}

Following Fatt's seminal work \cite{fatt1956network,xiong2016review}, C-PNM has become an efficient tool for rapid simulation of porous media flows.
C-PNM models a porous medium by using a graph whose nodes represents \emph{pore bodies} while edges represents \emph{pore throats}.
% C-PNM represents a porous medium by a graph whose nodes are \emph{pore bodies} and edges are \emph{pore throats}. 
Starting from a voxel image of the void space, the network is typically extracted using maximal ball \cite{dong2009pore}, medial axis \cite{lindquist1996medial,yi2017pore,zhang2024medial}, or watershed based methods \cite{gostick2017versatile}. The maximal ball approach, for example, identifies pore bodies as the largest inscribed spheres in the \emph{void} phase (red circles in Figure \ref{fig:model_pb}(b)) and connects neighboring spheres through throats whose medial axes are shown as green lines in Figure \ref{fig:model_pb}(b). Each adjacent pair $(i,j)$ thus produces a throat with effective length $L_{ij}$ and a cross-sectional descriptor (e.g., radius $r_{ij}$ for the 3-D circular pipe model, and channel width $w_{ij}$ for the 2-D planar channel model). 

For steady incompressible single-phase flow at low Reynolds number, it is often observed that the viscous force within throats are typically significantly larger than that within pore bodies.  As a result, C-PNM assumes that pressure drops concentrate in throats while pressure does not vary within pore bodies, as a reasonable approximation.
%For steady incompressible single-phase flow at low Reynolds number, C-PNM assumes that pressure drops concentrate in throats. 
The volumetric flux in throat $(i,j)$ is modeled by a one-dimensional Hagen-Poiseuille law
\begin{equation}
q^{\mathrm{PNM}}_{ij} = -\,\dfrac{\pi r_{ij}^{4}}{8\,\mu\,L_{ij}}\big(p^{\mathrm{PNM}}_{i}-p^{\mathrm{PNM}}_{j}\big),
    \qquad (i,j)\in \mathcal E,
\label{eq:hp_3d}
\end{equation}
where $p_i^{\rm PNM}$ is the pressure at pore $i$, $\mu$ is the dynamic viscosity, and $\mathcal{E}$ denotes the node set of the pore graph. We remark that equation \eqref{eq:hp_3d} formally describes one-dimensional flow, but is derived from 3-D circular pipe flow. Mass conservation at every node $i$ produces a sparse linear system
$\sum_{j} q_{ij}^{\mathrm{PNM}}=0$ that can be solved efficiently by using direct or iterative solvers. 

In this work, we consider only 2-D test cases; 
therefore, whenever a 1-D closure is used, all numerical simulations employ the following Hagen–Poiseuille relation derived from  
2-D planar channel flow with channel width $w_{ij}$
\begin{equation*}
q^{\mathrm{PNM}}_{ij}
= -\,\frac{w_{ij}^{3}}{12\,\mu\,L_{ij}}\Big(p^{\mathrm{PNM}}_{i}-p^{\mathrm{PNM}}_{j}\Big), \qquad (i,j)\in \mathcal{E}.
\label{eq:hp_2d}
\end{equation*}

%\begin{equation*}
%q_{ij}^{\mathrm{PNM}}
%  = -\frac{\pi r_{ij}^{4}}{8\mu\,L_{ij}}\,
%    \bigl(p_i^{\mathrm{PNM}}-p_j^{\mathrm{PNM}}\bigr),
%    \qquad (i,j)\in \mathcal E,
%\end{equation*}

Although C-PNM is computationally attractive,
its accuracy hinges on idealized cross-sections and 1-D closures. 
% its simplifying assumptions (C-PNM limitations ??) ... 
%is known to lose accuracy in
%tortuous or highly constricted throats where the circular-tube
%assumption of \eqref{eq:hp} is violated. 
 The remainder of this section
develops a hybrid strategy that retains FE fidelity without incurring
global FE costs.

\subsection{Reference finite-element formulation}\label{sec:fem}
We now formulate the reference finite-element (FE) solution of \eqref{eq:stokes-global}. This full FE solution  provides (i) a high-fidelity baseline for validation (Section \ref{sec:num_results}) of DD-PNM (Section \ref{sec:ddpnm}) and (ii) the local operators required to assemble the global interface system in DD-PNM.

 We define the velocity and pressure spaces
\begin{equation*}
\mathcal V = \{\mathbf{v}\in [H^1(\Omega)]^2: \mathbf{v} = \mathbf{0} \text{ on } \Gamma_{\rm dir}\},\quad \mathcal Q:=L^2(\Omega).  
\end{equation*}
 We introduce the weak formulation: 
find $(\mathbf u,p)\in \mathcal V\times \mathcal Q$ such that
\begin{equation}
\left\{
\begin{array}{ll}
  a(\mathbf u,\mathbf v) + b(\mathbf v,p) = g(\mathbf v)
  &\forall\,\mathbf v\in \mathcal V,\\
  b(\mathbf u,q)                         = 0
  &\forall\,q\in \mathcal Q,
  \end{array}
  \right.
  \label{eq:weak_stokes_glob}
\end{equation}
where the bilinear forms
$$a(\mathbf u,\mathbf v)=\nu\int_\Omega\nabla\mathbf u\!:\!\nabla\mathbf v\,dx,\quad
b(\mathbf v,q)= -\int_\Omega q\,\nabla\!\cdot\!\mathbf v\,dx,$$
and the linear functional 
$$g(\mathbf v)=\int_{\Gamma_{\rm{in}}\cup \Gamma_{\rm{out}}}
\boldsymbol{\sigma}(\mathbf u,p)\mathbf{n}\cdot \mathbf{v}\,dx.
$$
We remark that, since $\mathbf{u} \in [H^1(\Omega)]^2$ and $p\in L^2(\Omega)$, the stress has no trace; hence the linear functional is used only for derivation purposes and does not appear in the final weak form.  
On the left and right boundaries $\Gamma_{\rm{in}}\cup \Gamma_{\rm{out}}$ we prescribe
\begin{equation}
\mathbf n\cdot\boldsymbol{\sigma}(\mathbf u,p)\mathbf{n} = -p_{\rm b},\quad \mathbf t\cdot\boldsymbol{\sigma}(\mathbf u,p)\mathbf{n} = 0\; (\forall\,\mathbf t\perp\mathbf n),
\label{eq:assumption_interface_stress}
\end{equation}
with $p_{\rm b}=p_{\rm in}$ on $\Gamma_{\rm in}$, $p_{\rm b}=p_{\rm out}$ on $\Gamma_{\rm out}$.
Equivalently, $\boldsymbol{\sigma}(\mathbf u,p)\mathbf{n} = -p_{\rm b}\mathbf{n}$, i.e., 
 the normal traction equals the specified boundary pressure and the tangential components of the traction vanish.  
 Substituting  \eqref{eq:assumption_interface_stress} into $g(\mathbf v)$  yields the explicit load
\begin{equation*}
g(\mathbf v)=-p_{\rm in}\int_{\Gamma_{\rm{in}}}
\mathbf{n}\cdot \mathbf{v}\,dx -p_{\rm out}\int_{\Gamma_{\rm{out}}}
\mathbf{n}\cdot \mathbf{v}\,dx.
\end{equation*}

We proceed to discretize the weak formulation \eqref{eq:weak_stokes_glob}. 
We employ the inf-sup stable Taylor-Hood pair
$\mathcal V_h = [\mathbb P_{2}]^2 \subset \mathcal V$, $\mathcal Q_h = \mathbb P_{1} \subset \mathcal Q$
on a body-fitted triangular mesh (Figure \ref{fig:model_pb}(d)). Let $\{\boldsymbol{\phi}_j\}_{j=1}^{n_{\rm u}}$ and $\{\psi_j\}_{j=1}^{n_{\rm p}}$ denote the velocity and pressure Lagrange bases, the FE approximations expand as
$$
\mathbf{u}_{h} = \sum\limits_{j=1}^{n_{\rm u}}U_j\boldsymbol{\phi}_j,\quad 
p_h = \sum\limits_{l=1}^{n_{\rm p}}P_l{\psi}_l,
$$
with coefficient vectors $\mathbf{U}=(U_1,\cdots, U_{n_{\rm u}})^{\top}$ and 
$\mathbf{P}=(P_1,\cdots, P_{n_{\rm p}})^{\top}$. The FE counterpart of \eqref{eq:weak_stokes_glob} in matrix form reads
\[
\begin{bmatrix}
A &  B^\top\\
 B &  0
\end{bmatrix}
\begin{bmatrix}
\mathbf U\\ \mathbf P
\end{bmatrix}
=
\begin{bmatrix}
\mathbf g\\ \mathbf 0
\end{bmatrix},
\qquad
 A_{jk}=a(\boldsymbol\phi_k,\boldsymbol\phi_j),\;
 B_{lj}=b(\boldsymbol\phi_j,\psi_l),\;\mathbf{g}_j = g(\boldsymbol\phi_j).
\]

While the global FE formulation delivers high fidelity, it becomes impractical on
 industry-scale images due to the prohibitive number of degrees of freedom. 
To retain accuracy at a reduced cost, we therefore decompose the domain
into independent pore-scale subdomains (Figure \ref{fig:model_pb}(c))
and replace the single global solver  with independent subdomain solvers, which we couple by assembling the local Dirichlet-to-Neumann maps $G_i$ 
 into the interface Schur complement $S$. 
This leads to the domain-decomposition pore-network method (DD-PNM) introduced in Section \ref{sec:ddpnm}.

\subsection{Domain-decomposition PNM (DD-PNM)}

\label{sec:ddpnm}

We partition $\Omega$ into $N$ non-overlapping Lipschitz subdomains $\Omega_{i},\;i=1,2,\cdots,N$. 
 For any pair $(i,j)$ of adjacent pores, the shared interface is $\Gamma_{ij}=\partial\Omega_{i}\cap\partial\Omega_{j}$, 
which is constructed by connecting the center of the interface to the nearest point on each of its two neigboring solid balls, while the center of the interface is located by identifying a saddle point of the distance function  
% which we locate by identifying saddle points of the distance function
  \cite{liu2024new}
 $$
 d(\mathbf x)=\min_{\mathbf y\in \rm{solid}}\Vert \mathbf x - \mathbf y \Vert,
 $$
 i.e., the distance from a void-space point $\mathbf x$ to the solid phase. In the geometry of Figure \ref{fig:model_pb}(c), this choice places the interfaces along the lines connecting neighboring solid centers (blue lines). Other interface-placement strategies are possible; a systematic study is beyond the scope of the present work. 
In Section \ref{sec:num_results} we
demonstrate robustness of the method with respect to interface placement: even when the
locations of $\Gamma_{ij}$ are perturbed, the computed results remain accurate.

We state the \texttt{main modeling assumption} underlying DD-PNM: on each internal interface $\Gamma_{\alpha}$, the traction is fully determined by a single scalar interface variable $p_{\alpha}$. Specifically, we assume 
\begin{equation}
\mathbf n\cdot\boldsymbol{\sigma}(\mathbf u,p)\mathbf{n} = -p_{\alpha},\quad \mathbf t\cdot\boldsymbol{\sigma}(\mathbf u,p)\mathbf{n} = 0\; (\forall\,\mathbf t\perp\mathbf n)\;\text{on}\;\Gamma_{\alpha},
\label{eq:assumption_interface_stress2}
\end{equation}
which enforces a uniform normal traction and zero tangential traction across each interface. This assumption is consistent with the inlet/outlet boundary treatment introduced in \eqref{eq:assumption_interface_stress}, and ensures that subdomain coupling is fully described by the set of interface variables $\{p_{\alpha}\}_{\alpha=1}^m$, $m$ denotes the number of internal faces. 

\begin{remark}
\label{rem:interface_alt}
An alternative to the traction-based interface condition \eqref{eq:assumption_interface_stress2} is to assume that the pressure is constant along each internal face, i.e., $p = p_\alpha$ on $\Gamma_\alpha$, and that the tangential component of the velocity vanishes:
$$
p = p_\alpha,\qquad \mathbf{u} \cdot \mathbf{t} = 0 \quad (\forall\, \mathbf{t} \perp \mathbf{n}) \quad \text{on }\; \Gamma_\alpha.
$$
This formulation also leads to a single scalar degree of freedom per interface and reflects the physical assumption of dominant normal flow between pores. A detailed comparison of different interface modeling assumptions will be the subject of future work.
\end{remark}

%  The present choice is adequate, here the pressure varies slightly, conform to our main assumption in the DD-PNM, \texttt{interface pressure is supposed to be constant}. 
  
%  Let  
%$\Gamma_\alpha,\;\alpha=1,\dots,m$,  
%denote the $m$ internal interfaces, and let  
%$p_\alpha$ be the (unknown) constant pressure prescribed on
%$\Gamma_\alpha$.
For each subdomain $\Omega_i$ we list its $m_i$ neighbouring
faces as
$\Gamma_{i,1},\dots,\Gamma_{i,m_i}$ and 
split them into  an 
        unknown-traction set $\mathcal U_i$ and a known-traction (inlet/outlet) set $\mathcal K_i$
\[
\mathcal U_i = \{\Gamma_{i,1},\dots,\Gamma_{i,m_i^{\rm u}}\},\qquad
\mathcal K_i = \{\Gamma_{i,m_i^{\rm u}+1},\dots,\Gamma_{i,m_i}\},
   \quad m_i^{\rm u}+m_i^{\rm k}=m_i.
\]
The associated local traction vectors are $
\mathbf p_i^{\rm u} =
\begin{bmatrix}
p_{i,1},\cdots, p_{i,m_i^{\rm u}}
\end{bmatrix}^{\top}$,\\
$\mathbf p_i^{\rm k} =
\begin{bmatrix}
p_{i,m_i^{\rm u}+1},\cdots, p_{i,m_i}
\end{bmatrix}^{\top}$,
$\mathbf p_i=\begin{bmatrix} p_{i,1},\cdots, p_{i,m_i}\end{bmatrix}^{\top}.
$
We define the local-to-global index maps
\begin{itemize}\setlength{\itemsep}{2pt}
\item
$\beta_i^{\rm u}:\{1,\dots,m_i^{\rm u}\}\rightarrow\{1,\dots,m\}$  
maps the $r$-th \emph{internal} face of $\Omega_i$ to its global
internal index:
$ p_{i,r}=p_{\beta_i^{\rm u}(r)}\; (r\le m_i^{\rm u})$.
\item
$\beta_i^{\rm k}:\{1,\dots,m_i^{\rm k}\}\rightarrow\{1,\dots,m^{\rm k}\}$  
maps the $s$-th \emph{boundary} face of $\Omega_i$
($s=m_i^{\rm u}+1,\dots,m_i$) to the global inlet/outlet list, whose
entries are prescribed and collected in
$\mathbf p^{\rm k}\in\mathbb R^{m^{\rm k}}$.
\end{itemize}
We define the Boolean restriction matrices
\[
(R_i^{\rm u})_{r,\beta_i^{\rm u}(r)}=1,\quad
R_i^{\rm u}\in\{0,1\}^{\,m_i^{\rm u}\times m},
\qquad
(R_i^{\rm k})_{s,\beta_i^{\rm k}(s)}=1,\quad
R_i^{\rm k}\in\{0,1\}^{\,m_i^{\rm k}\times m^{\rm k}}.
\]
With the global traction vectors
$
\mathbf p=[p_1,\dots,p_m]^{\!\top}
$ and $\mathbf p^{\rm k}$
we obtain
\[
\mathbf p_i^{\rm u}=R_i^{\rm u}\,\mathbf p,\quad \mathbf p_i^{\rm k}=R_i^{\rm k}\,\mathbf p^{\rm k}. 
\] 
These matrices and vectors will be subsequently used to
assemble the global Schur complement and its forcing term.

\paragraph{Local Stokes solve and Dirichlet-to-Neumann map}
Within each subdomain $\Omega_i$  we again solve the steady Stokes equations
\begin{equation}
\left\{
\begin{array}{ll}
-\nu\Delta\mathbf{u}_i + \nabla \pi_i = \mathbf{0}, &\text{in } \Omega_i, \\
\nabla\!\cdot\!\mathbf{u}_i = 0, &
\text{in } \Omega_i, \\
\mathbf u_i = \mathbf 0,
&\text{on } \partial\Omega_i^{\text{wall}}, \\
\boldsymbol{\sigma}(\mathbf{u}_i,\pi_i)\mathbf{n}_{i,r} = -p_{i,r} \mathbf{n}_{i,r}, & \text{on } \Gamma_{i,r},\; r=1,\dots,m_i.
\end{array}
\right.
\label{eq:stokes-local}
\end{equation}
Here, the symbol $\pi$ is used to represent the local  pressure field to avoid confusion with the interface normal traction $p$.

%where $\mathcal J_i$ is the set of neighbours of $i$. 

Using the Taylor-Hood discretization of 
Section \ref{sec:fem}, we obtain the linear system 
%  Prescribing the $m_i$ face pressures yields
\begin{equation}
\begin{bmatrix}
A_i & B_i^{\!\top}\\[2pt]
B_i & 0
\end{bmatrix}
\begin{bmatrix}\mathbf U_i\\[2pt]\mathbf P_i\end{bmatrix}
=
\sum_{r=1}^{m_i}p_{i,r}
\begin{bmatrix}\mathbf f_i^{(r)}\\[2pt]0\end{bmatrix},
\label{eq:stokes_local_matrix}
\end{equation}
where $\mathbf{U}_i$, $\mathbf{P}_i$ are the coefficient vectors of $\mathbf{u}_i^h$, $\pi_i^h$, and
$\mathbf f_i^{(r)}\in\mathbb{R}^{n_{i,\rm u}}$ is the Neumann load vector for a unit normal traction on the interface $\Gamma_{i,r}$, with components
  \begin{equation}
  \bigl(\mathbf f_i^{(r)}\bigr)_j
    \;=\;
    -\int_{\Gamma_{i,r}}
      \boldsymbol{\phi}_{i,j}\cdot\mathbf{n}_{i,r}\,dx.
      \label{eq:neumann_load}
  \end{equation}   
Because the right-hand side is linear in   $p_{i,r}$, the discrete velocity depends linearly on those tractions
\begin{equation}
\mathbf{u}_i^h = \sum_{r=1}^{m_i} p_{i,r}\, \widetilde{f}_i^{(r)},
\qquad
\mathbf{U}_i = \sum_{r=1}^{m_i} p_{i,r}\, \widetilde{\mathbf{f}}_i^{(r)},
\label{eq:vel_pressure}
\end{equation}
where $\widetilde{f}_i^{(r)}$ is the velocity response function to a unit normal traction on $\Gamma_{i,r}$ and zero traction on all other faces,
and $\widetilde{\mathbf{f}}_i^{(r)} \in \mathbb{R}^{n_{i,\rm u}}$ its coefficient vector (its explicit expression is given in Lemma \ref{lem:disc_incomp_ftilde}).

\begin{remark}
Although an explicit formula for $\widetilde{\mathbf f}_i^{(r)}$ is provided in Lemma \ref{lem:disc_incomp_ftilde}, in practice we obtain $\widetilde{\mathbf f}_i^{(r)}$ by solving \eqref{eq:stokes_local_matrix} with unit normal traction on $\Gamma_{i,r}$ and zero on the remaining interfaces, and then extracting the
velocity field from the resulting solution.
\end{remark}

The outward flux through face $\Gamma_{i,k}$ is defined by 
\[Q_{i,k} = \int_{\Gamma_{i,k}} \mathbf u_i^h \cdot \mathbf{n}_{i,k}\,dx.\]  
Substituting the velocity expansion \eqref{eq:vel_pressure} and invoking the 
 Neumann load vector \eqref{eq:neumann_load} gives
\begin{equation}
Q_{i,k} = -\,[\mathbf{f}_i^{(k)}]^{\top} \mathbf{U}_i =  \sum_{r=1}^{m_i} p_{i,r}\,g_{i,kr},\quad g_{i,kr}:= - [\mathbf f_i^{(k)}]^{\top}\,\widetilde{\mathbf f}_i^{(r)}.
\label{eq:flux_and_g}
\end{equation}
%    Introducing the row vector $C_{i,k} = - [f_i^{(k)}]^{\top}$, we obtain
%\begin{equation}
%Q_{i,k} = C_{i,k}\mathbf{u}_i^h = \sum_{r=1}^{m_i} p_{i,r}\,g_{i,kr},\quad g_{i,kr}:= C_{i,k}\,\widetilde{f}_i^{(r)}.
%\label{eq:flux_and_g}
%\end{equation}
Hence $g_{i,kr}$ is the normal flux across $\Gamma_{i,k}$ produced by a unit normal traction applied on $\Gamma_{i,r}$ while all other interface tractions are set to zero. 
Collecting $g_{i,kr}$ gives the local (discrete) \emph{Dirichlet-to-Neumann} (DtN), or
\emph{Steklov-Poincaré}, operator
% the symmetric positive-definite (which will be proved in Theorem ??) matrix
\begin{equation}
  G_i
  \;=\;
  \bigl[g_{i,kr}\bigr]_{k,r=1}^{m_i}
  \;\in\;\mathbb{R}^{m_i\times m_i},
  \label{eq:dtn_local}
\end{equation}
which  maps an \emph{interface-traction vector} to the
\emph{flux vector}. 

Similar to local traction vectors, we define the local flux vectors
$$
\mathbf Q_i^{\rm u}=[Q_{i,1},\ldots,Q_{i,m_i^{\rm u}}]^{\top},\quad \mathbf Q_i^{\rm k}=[Q_{i,m_i^{\rm u}+1},\ldots,Q_{i,m_i}]^{\top}, \quad \mathbf Q_i=[Q_{i,1},\ldots,Q_{i,m_i}]^{\top}.
$$
The local traction vector and local flux vector can thus be related by
\begin{equation}
\mathbf Q_i = G_i\,\mathbf p_i. 
\label{eq:relation_Qi_pi}
\end{equation}
Splitting into known/unknown blocs yields
\[
\begin{bmatrix}
\mathbf Q_i^{\rm u}\\[2pt]\mathbf Q_i^{\rm k}
\end{bmatrix}
=
\begin{bmatrix}
G_i^{\rm uu} & G_i^{\rm uk}\\
G_i^{\rm ku} & G_i^{\rm kk}
\end{bmatrix}
\begin{bmatrix}
\mathbf p_i^{\rm u}\\[2pt]\mathbf p_i^{\rm k}
\end{bmatrix},
\]
and the internal face flux can be computed as
$
\mathbf Q_i^{\rm u} = G_i^{\rm uu}\, \mathbf p_i^{\rm u}  +  G_i^{\rm uk}\,\mathbf p_i^{\rm k}.
$

\paragraph{Global coupling and Schur complement}

Using the restriction matrices, the subdomain flux relation can be lifted to global indices
 $$
 (R_i^{\rm u})^{\!\top}\mathbf Q_i^{\rm u} = (R_i^{\rm u})^{\!\top}G_i^{\rm uu}\,\mathbf p_i^{\rm u} + (R_i^{\rm u})^{\!\top}G_i^{\rm uk}\,\mathbf p_i^{\rm k} = (R_i^{\rm u})^{\!\top}G_i^{\rm uu}\,R_i^{\rm u} \mathbf p + (R_i^{\rm u})^{\!\top}G_i^{\rm uk}\,R_i^{\rm k} \mathbf p^{\rm k}.
 $$ 
Summing over all subdomains and enforcing mass balance on every interface 
 gives the global Schur complement system
\begin{equation}
0
  = -\sum_{i}^N(R_i^{\rm u})^{\!\top}\mathbf Q_i^{\rm u}
  = \underbrace{-\sum_{i}^N (R_i^{\rm u})^{\!\top}G_i^{\rm uu}\,R_i^{\rm u} }_S \mathbf p - \underbrace{\,\sum_{i}^N (R_i^{\rm u})^{\!\top}G_i^{\rm uk}\,R_i^{\rm k} \mathbf p^{\rm k}}_{\mathbf F},
  \label{eq:global_schur}
\end{equation}
or simply
\[
S\,\mathbf p=\mathbf F.
\]
where the forcing vector $\mathbf F$ collects the contributions from faces with prescribed (inlet/outlet) tractions.  

This construction completes the DD-PNM formulation: local
Stokes problems yield DtN blocks $G_i$, whose global assembly
\eqref{eq:global_schur} enforces flux continuity while keeping the
system size proportional to the number of internal faces.

\begin{remark}
With $g_{i,kr}=-[\mathbf f_i^{(k)}]^{\!\top}\widetilde{\mathbf f}_i^{(r)}$,
we have the energy identity
$\mathbf p_i^{\!\top}G_i\mathbf p_i
= -\,\mathbf U_i^{\top}A_i\mathbf U_i\le 0$,
so $G_i$ is symmetric negative semi-definite with the uniform-traction vector
in its kernel (cf.\ Section \ref{sec:properties_DD-PNM}). We therefore flip the
sign in the global assembly to obtain a symmetric positive-definite (SPD) Schur complement.

\end{remark}

%  \item $A_i$ – the local stiffness matrix obtained by discretising the inertia-free
%        Stokes equations on $\Omega_i$; $A_i$ is symmetric positive-definite (SPD).
%  \item $F_\alpha^{(i)}$ – the interface load vector on $\Omega_i$ that corresponds to
%        applying a \emph{unit} pressure on $\Gamma_\alpha$
%        (all other interfaces being set to zero).
%  \item $u_i^{(\alpha)} = A_i^{-1}F_\alpha^{(i)}$ – the velocity response in $\Omega_i$
%        produced by a unit pressure on $\Gamma_\alpha$.

\paragraph{Solution algorithm}
Algorithm \ref{alg:ddpnm} summarizes the DD-PNM workflow. All subdomain
computations are independent and executed in parallel. Local responses are
obtained by solving the discrete Stokes problem on each
subdomain.

\begin{algorithm}[h]
\caption{DD-PNM solver (assemble - couple - recover)}
\label{alg:ddpnm}
\begin{algorithmic}[1]
\STATE \textbf{Input:} subdomain meshes and FE pairs $(V_i^h,Q_i^h)$; prescribed tractions $\mathbf p^{\rm k}$.
\STATE \textbf{Output:} unknown interface tractions $\mathbf p$, interface fluxes, and (optionally) local fields.
\vspace{0.25em}
\FORALL{pore bodies $\Omega_i$ \textbf{in parallel}}
  \STATE Assemble local operators $A_i$, $B_i$ and unit normal traction load vectors $\{\mathbf f_i^{(r)}\}_{r=1}^{m_i}$.
  \FORALL{local faces $r=1,\dots,m_i$}
    \STATE Solve the local Stokes system for the unit normal traction  response:
           \begin{equation}
             \begin{bmatrix} A_i & B_i^{\top}\\[2pt] B_i & 0 \end{bmatrix}
             \begin{bmatrix} \widetilde{\mathbf f}_i^{(r)} \\[2pt] \mathbf X_i^{(r)} \end{bmatrix}
             =
             \begin{bmatrix} \mathbf f_i^{(r)} \\[2pt] 0 \end{bmatrix}.
             \label{eq:sys_unit_pressure_response}
           \end{equation}
           (Here $\widetilde{\mathbf f}_i^{(r)}$ is the velocity component, and $\mathbf X_i^{(r)}$ is the pressure component.)
  \ENDFOR
  \STATE Form the local DtN blocks $G_i$ with entries 
         $g_{i,kr} = -\,[\mathbf f_i^{(k)}]^{\top}\,\widetilde{\mathbf f}_i^{(r)}$ and split $G_i$
         into blocks $G_i^{\rm uu}$  and $G_i^{\rm uk}$.
  \STATE Accumulate global coupling with Boolean maps $R_i^{\rm u},R_i^{\rm k}$:
         \[
           S \;{-}{=}\; (R_i^{\rm u})^{\!\top} G_i^{\rm uu} R_i^{\rm u},
           \qquad
           \mathbf F \;{+}{=}\; (R_i^{\rm u})^{\!\top} G_i^{\rm uk} R_i^{\rm k}\,\mathbf p^{\rm k}.
         \]
\ENDFOR
\STATE Solve the global interface system \quad $S\,\mathbf p = \mathbf F$.
\STATE {(Optional)} Reconstruct local fields (in parallel):
       \[
         \mathbf U_i
           = \sum_{r=1}^{m_i} p_{i,r}\,\widetilde{\mathbf f}_i^{(r)},
\qquad             
   \mathbf P_i
           = \sum_{r=1}^{m_i} p_{i,r}\,\mathbf X_i^{(r)},           
             \qquad
         \mathbf Q_i^{\rm u} = G_i^{\rm uu}\mathbf p_i^{\rm u} + G_i^{\rm uk}\mathbf p_i^{\rm k}.
       \]
\end{algorithmic}
\end{algorithm}

%\begin{algorithm}
%\caption{Build tree}
%\label{alg:buildtree}
%\begin{algorithmic}
%\STATE{Define $P:=T:=\{ \{1\},\ldots,\{d\}$\}}
%\WHILE{$\#P > 1$}
%\STATE{Choose $C^\prime\in\mathcal{C}_p(P)$ with $C^\prime := \operatorname{argmin}_{C\in\mathcal{C}_p(P)} \varrho(C)$}
%\STATE{Find an optimal partition tree $T_{C^\prime}$ }
%\STATE{Update $P := (P{\setminus} C^\prime) \cup \{ \bigcup_{t\in C^\prime} t \}$}
%\STATE{Update $T := T \cup \{ \bigcup_{t\in\tau} t : \tau\in T_{C^\prime}{\setminus} \mathcal{L}(T_{C^\prime})\}$}
%\ENDWHILE
%\RETURN $T$
%\end{algorithmic}
%\end{algorithm}

\paragraph{Computational cost}
Each subdomain $\Omega_i$ has a small number of interfaces $m_i=\mathcal O(1)$, so the local DtN matrix $G_i\in\mathbb R^{m_i\times m_i}$ is tiny and the assembled global matrix $S$ is extremely sparse. Building $G_i$ requires $m_i$ \emph{independent} local Stokes solves; this work is embarrassingly parallel across subdomains. Assembly of $S$ and $\mathbf F$ is linear in $N$ (the number of pore-subdomains), and the size of the global system equals the number of unknown interface tractions rather than the number of FE nodes.

\begin{remark}
Once $\mathbf p$ is computed, interface fluxes $\mathbf Q_i^u$ follow directly from the DtN maps. When detailed fields are
needed, fine-scale velocity and pressure can be recovered locally from the
precomputed unit normal traction responses $\widetilde{\mathbf f}_i^{(r)}$ and $\mathbf X_i^{(r)}$ without any
additional global solve (see Algorithm \ref{alg:ddpnm}).

\end{remark}

\subsection{Recovering C‑PNM from the DD‑PNM DtN maps}

\label{sec:recovering_CPNM}

We now show how to recover a classical pore–network model (C–PNM) directly from the DD–PNM Dirichlet–to–Neumann (DtN) maps computed in Section~\ref{sec:ddpnm}, thereby providing a principled bridge between mesh-based and graph-based models within a unified algebraic framework. 
C-PNM assigns a scalar conductivity to each throat (typically from geometric surrogates) and writes face fluxes in terms of the difference between neighboring pore–body pressures and a pore–body pressure at the current node. 
In our setting, by contrast, we \emph{calibrate} the conductivities associated with the faces of each pore directly from the FE subdomain responses, so that the resulting linear face–flux map reproduces the DD-PNM interface map.
 This calibration preserves the FE-resolved intra-pore physics while remaining fully compatible with established C–PNM workflows.
 
\begin{figure}[htbp!]
		\centering
		\includegraphics[scale = 0.36]{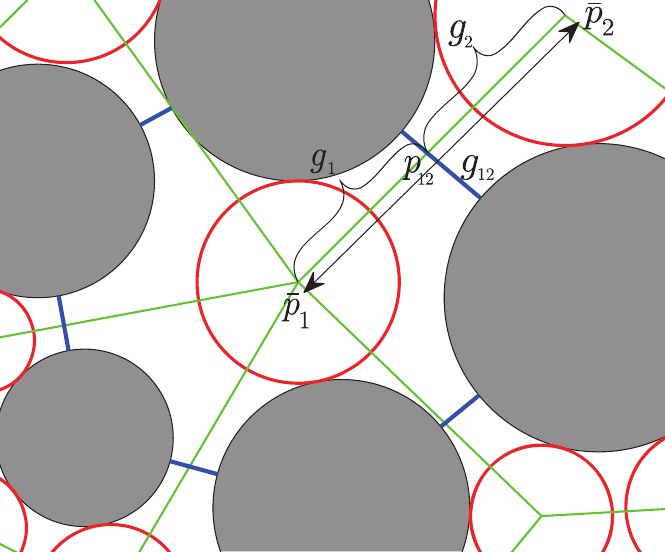}
		\caption{C-PNM versus the calibration of C-PNM from DD-PNM of a throat.}
		\label{fig:PNM_model_marked}
	\end{figure}

As an illustrative example, Figure \ref{fig:PNM_model_marked} provides a zoomed-in view of two neighboring pores previously shown in Figure \ref{fig:model_pb}(b). In the C-PNM formulation, 
the flux through the throat is $\widehat q_{12}=g_{12}(\bar p_1-\bar p_2)$ in terms of the \emph{pore-body pressures} $(\bar p_1,\bar p_2)$ and a single throat conductivity, denoted by $g_{12}$ in the figure.  
In contrast, our DD-PNM to C-PNM calibration introduces an \emph{interface scalar} $p_{12}$ on the shared interface, which in our formulation is the normal traction (unit: Pa). 
This divides the throat into two segments, referred as \emph{half-throats}, each connecting the shared interface to one of the pore bodies. 
We compute the two half-throat conductivities $(g_1,g_2)$ from the FE subdomain responses so that the face–flux relations 
$\widehat q_1=g_1(\bar p_1-p_{12})$ and $\widehat q_2=g_2(p_{12}-\bar p_2)$ reproduce the DD-PNM Dirichlet-to–Neumann map. 
The half-throat conductivities are then combined by the harmonic average to recover the single $g_{12}$ used in C-PNM, as detailed  later in this section. 

For a pore subdomain $\Omega_i$ with $m_i$ faces, introduce the (nonnegative) face-wise conductivities and the diagonal matrix 
$$D_i:=\mathrm{diag}(g_{i,1},\ldots,g_{i,m_i}),\quad g_{i,r}\geq 0\; (r=1,\cdots,m_i).$$ 
Let $\mathbf{p}_i = [\,p_{i,1},\dots,p_{i,m_i}\,]^\top$ denotes the vector of interface tractions (in Pa) of $\Omega_i$. The C-PNM computes the vector of face fluxes $\mathbf{q}_i$ as
\begin{equation}
\mathbf{q}_i=D_i(\mathbf p_i-\bar{p}_i\,\mathbf{1}),
\label{eq:cpnm_face_flux}
\end{equation}
where $\mathbf{1}\in\mathbb{R}^{m_i}$ is the all-ones vector
 and $\bar{p}_i$ is the node (pore-body)  pressure used by the recovered network. Relation \eqref{eq:cpnm_face_flux} states that the flux on face $r$ is 
$q_{i,r}=g_{i,r}(p_{i,r}-\bar{p}_i)$, where 
 the difference $p_{i,r}-\bar{p}_{i}$ corresponds to the effective normal-traction drop across the half-throat segment
 % (here referred as half-throat for brevity) 
 connecting the  pore-body to the interface, and $g_{i,r}$ denotes the fitted conductivity. 
 % the vector of interface scalar variables on the faces of
%$\Omega_i$ (in our formulation these are normal tractions, in Pa).
%where 
%$\bar{p}_i$ is the (unknown) pore-body pressure, $\mathbf{1}$ denotes the $m_i$-dimensional all-ones vector. We remark that  C-PNM takes $p_i$ as neighboring pore-body pressures, and $g_{i,r}$ as throat conductivities; however, in our DD-PNM setting, ${p}_i$ contains interface pressures, and $g_{i,r}$ denotes the conductivities associated with. 
%Here, $p_i$ denotes face-wise interface tractions (Pa) and  
%Together with , this difference reproduces the FE-computed face fluxes through the relation 
Mass balance $\mathbf{1}^{\top}\mathbf q_i=0$ gives
$$
\bar{p}_i=\frac{\mathbf{1}^\top D_i\mathbf p_i}{\mathbf{1}^\top D_i\mathbf{1}}.
$$
Eliminating $\bar{p}_i$ from the C-PNM relation gives a linear map on face tractions,
\begin{equation}
\mathbf q_{i}=G_{i}^{\star}(D_{i})\mathbf p_{i},\quad G_{i}^{\star}(D_{i})=D_{i}-\frac{D_{i}\mathbf{1}\mathbf{1}^{\top}D_{i}}{\mathbf{1}^{\top}D_{i}\mathbf{1}}.
\label{eq:Gi_star}
\end{equation}
%Entry-wise form:
%$$
%\left[G_i^\star(D_i)\right]_{\omega r}=
%\left\{
%\begin{array}{ll}
%-\frac{g_{i,\omega}g_{i,r}}{\sum_l g_{i,l}},&\omega\neq r,\\
%g_{i,\omega}-\frac{g_{i,\omega}^2}{\sum_l g_{i,l}},&\omega=r.
%\end{array}
%\right.
%$$
On the other hand, DD-PNM provides the FE-accurate interface map
 $\mathbf Q_i=G_i\mathbf p_i$ for each pore subdomain (see equation \eqref{eq:relation_Qi_pi}).  
 We therefore determine $D_i$ so that $G_i^\star(D_i)$ \emph{matches} $G_i$ (up to the sign convention induced by outward normals)
\begin{equation}
G_i \approx -G_i^*(D_i), 
\end{equation}
ensuring that $\mathbf q_i \approx \mathbf{Q}_i$ after calibration. 
Concretely, we compute the half-throat conductivities by solving the nonnegative nonlinear least-squares problem
$$\min_{g_{i,r}\geq0}\left.\left\|G_i+G_i^\star(D_i)\right\|_F^2,\right.$$
where $\left\|\cdot \right\|_F$ denotes the Frobenius norm. 
Since both $G_i$ and $G_i^\star$ have zero row sum (the latter is immediate from \eqref{eq:Gi_star}; for $G_i$ see Lemma \ref{lem:const_kernel}), the diagonal entries are determined by the off-diagonals. It therefore suffices to fit only the off-diagonal entries: 
\begin{equation}\label{eq:off-fit}
\min_{g_{i,r}\ge 0}\;\;
\big\|\, \mathrm{Off}(G_i)\;+\;\mathrm{Off}\!\big(G_i^\star(D_i)\big)\,\big\|_{F}^2,
\qquad 
\mathrm{Off}(A):=A-\mathrm{diag}(A).
\end{equation}
Let $
T_i \;:=\; \mathrm{Off}(G_i) $
% \ge 0 $ 
and $\mathbf g_i=(g_{i,1},\cdots,g_{i,m_i})^{\top}$. 
\begin{remark}
 Intuitively, in a pore domain, if the pressure on a single face $r$ is increased while all other face pressures are held fixed, 
 the outward fluxes on the other faces  increase; that is,
$g_{i,kr} \geq 0 $ for $k\neq r$,  or equivalently $T_i\geq 0$. This off-diagonal nonnegativity is consistent with the DtN blocks computed in our test cases. 
\end{remark}
\noindent Introduce the scale-free unknown $\mathbf x_i:=\frac{\mathbf g_i}{\sqrt{\mathbf{1}^{\top}\mathbf g_i}}$ ($\mathbf g_i$ can be recovered from $\mathbf x_i$ by $\mathbf g_i=(\mathbf{1}^{\top}\mathbf x_{i})\mathbf x_i)$), so that
\[
\mathrm{Off}\!\big(G_i^\star(D_i)\big)
\;=\;
\mathrm{Off}\big(-\mathbf x_i \mathbf x_i^\top\big).
\]
Thus, the optimization problem that we aim to solve becomes
\begin{equation}\label{eq:masked-rank1}
\min_{\mathbf x\ge 0}\;\; 
f(\mathbf x)\;:=\;\big\|\, \mathbf x \mathbf x^\top - T_i \,\big\|_{F,\mathrm{off}}^2,
\qquad 
\|A\|_{F,\mathrm{off}}^2:=\sum_{r\neq s} A_{rs}^2.
\end{equation}

We solve \eqref{eq:masked-rank1} using the projected gradient descent (PGD) method \cite{beck2017first}. For an effective initializer, 
observe that \eqref{eq:masked-rank1} is rank-1 approximation problem for the symmetric nonnegative matrix $T_i$. 
Let $\mathbf v\geq0$ be the Perron eigenvector \cite{horn2012matrix} of $T_i$ and 
% see perron_frobenius.pdf
set $M_0:=\frac{\mathbf v\mathbf v^\top}{\mathbf{1}^\top \mathbf v}$.
The best scalar multiple of $M_0$ in the off-diagonal Frobenius norm is
$$\alpha^*\::=\:\arg\min_{\alpha\geq0}\:\left\|T-\alpha M_0\right\|_{F,\mathrm{off}}^2\:=\:\frac{\langle T,\:M_0\rangle_{\mathrm{off}}}{\|M_0\|_{F,\mathrm{off}}^2},\quad
\langle A,B\rangle_{\mathrm{off}}:=\sum_{i\neq j}A_{ij}B_{ij}.$$ 
We then set the initializer of the PGD method as $\mathbf x^{(0)}    = \sqrt{\alpha^*}\,\mathbf v$. 

If pores $i$ and $j$ share face $\omega$,  the corresponding throat conductivity is assembled from the two half-throat conductivities on that interface by the harmonic average
\[
g_\omega \;=\; \Big(g_{i,\omega}^{-1}+g_{j,\omega}^{-1}\Big)^{-1},
\]
and the C-PNM node equations (in the unknowns pore-body pressures  $\{\bar p_i\}$) follow from mass conservation. 
This DD-PNM to C-PNM reduction produces a \emph{calibrated} network representation: it preserves the FE-resolved intra-pore physics through the fitted half-throat conductivities and,
 once calibrated, can be used within standard C-PNM pipelines without further reference to the FE model.

\section{Mathematical properties of DD-PNM}
\label{sec:properties_DD-PNM}
In this section, we establish the solvability and mass-conservation properties of the DD-PNM formulation (see Theorem \ref{thm:ddpnm}). We begin by collecting several auxiliary lemmas that will be used in the proof.

\begin{lemma}[Local well-posedness]
      For every interface-traction vector \\
      $\mathbf p=(p_1,\dots,p_m)^{\top}$,
      the discrete Stokes subproblem in each $\Omega_i$ admits a unique solution. 
\label{lem:local_wellposed}
\end{lemma}

\begin{proof}
For each subdomain the saddle-point matrix
$
  K_i=\begin{bmatrix}A_i & B_i^{\top}\\ B_i & 0\end{bmatrix}
$ in~\eqref{eq:stokes_local_matrix} 
is invertible because
(i) $A_i$ is symmetric positive-definite (SPD) and
(ii) the pair $(V_i^h,Q_i^h)$ satisfies the discrete inf-sup condition.
Hence the local finite-element Stokes problem has a unique solution
$(\mathbf U_i,\mathbf P_i)$ for any given interface tractions $\mathbf p$.
\end{proof}

\begin{lemma}[Discrete incompressibility of the unit responses]
\label{lem:disc_incomp_ftilde}
For each local interface index $r$, the elementary velocity response $\widetilde{\mathbf f}_i^{(r)}$ (see Section \ref{sec:ddpnm}) is discretely divergence-free: 
$$
B_i\widetilde{\mathbf f}_i^{(r)} = 0.
$$
\end{lemma}
\begin{proof}
We first derive the explicit expression of $\widetilde{\mathbf f}_i^{(r)}$. From the first block row of \eqref{eq:stokes_local_matrix} we have
\begin{equation}
\mathbf{U}_i
=
A_i^{-1}
\Bigl(
\sum_{r=1}^{m_i} p_{i,r}\,\mathbf f_i^{(r)}
-
B_i^{\top}\mathbf P_i
\Bigr).
\label{eq:stokes_u}
\end{equation}
Imposing the discrete divergence-free constraint $B_i\mathbf{U}_i=0$ gives
\begin{equation}
\mathbf P_i
=
\sum_{r=1}^{m_i} p_{i,r}\;
M_i^{-1}B_iA_i^{-1}\mathbf f_i^{(r)},\quad M_i \;:=\; B_iA_i^{-1}B_i^{\top}\; \text{(SPD)}.
\label{eq:pi}
\end{equation}
Substituting \eqref{eq:pi} into equation \eqref{eq:stokes_u} yields the linear representation
\begin{align*}
\mathbf{U}_i
\;=\;
\sum_{r=1}^{m_i} p_{i,r}\,\widetilde{\mathbf f}_i^{(r)}, 
\end{align*}
with
\begin{align}
\widetilde{\mathbf f}_i^{(r)}
:=
A_i^{-1}\mathbf f_i^{(r)} - A_i^{-1}
B_i^{\top}
M_i^{-1}B_iA_i^{-1} \mathbf f_i^{(r)}.
\label{eq:f_tilde_def}
\end{align}
A direct calculation proves the claim: 
\begin{equation}
\begin{aligned}
B_i\widetilde{\mathbf f}_i^{(r)}
 = B_iA_i^{-1}\mathbf f_i^{(r)}
    -\underbrace{B_iA_i^{-1}B_i^{\!\top}}_{M_i}
     M_i^{-1}B_iA_i^{-1}\mathbf f_i^{(r)}
 = 0 .
\end{aligned}
\label{eq:f_tilde_incomp}
\end{equation}
\end{proof}

\begin{remark}
Lemma \ref{lem:disc_incomp_ftilde}
follows directly from \eqref{eq:sys_unit_pressure_response}: the second block row yields $B_i\,\widetilde{\mathbf f}_i^{(r)}=0$. However, here we present a detailed derivation to obtain the explicit expression of $\widetilde{\mathbf f}_i^{(r)}$, which will be used in the proof of Lemma \ref{lem:energy_identity}. 
\end{remark}

\begin{lemma}[Symmetry and energy identity for the local DtN map]
\label{lem:energy_identity}
Let $G_i$ be defined by $g_{i,kr}=-\,[\mathbf f_i^{(k)}]^{\top}\,\widetilde{\mathbf f}_i^{(r)}$.
Then $G_i$ is symmetric and, for any local traction vector
$\mathbf p_i=[p_{i,1},\dots,p_{i,m_i}]^{\top}$,
\begin{equation}
\mathbf p_i^{\!\top}G_i\,\mathbf p_i
  = -\,\mathbf U_i^{\!\top}A_i\,\mathbf U_i. 
\label{eq:local_energy_identity}
\end{equation}
In particular, $G_i$ is negative semi-definite.
\end{lemma}

\begin{proof}
By Lemma~\ref{lem:disc_incomp_ftilde}, $B_i\widetilde{\mathbf f}_i^{(r)}=0$. Using
\eqref{eq:f_tilde_def} and the symmetry of $A_i$,
\[
[\widetilde{\mathbf  f}_i^{(k)}]^{\top}A_i\widetilde{\mathbf f}_i^{(r)}
=\Bigl(\mathbf f_i^{(k)}-B_i^{\top}M_i^{-1}B_iA_i^{-1}\mathbf f_i^{(k)}\Bigr)^{\!\top}\widetilde{\mathbf f}_i^{(r)}
= [\mathbf f_i^{(k)}]^{\!\top}\widetilde{ \mathbf f}_i^{(r)}.
\]
Hence $g_{i,kr}=-\,[\mathbf f_i^{(k)}]^{\!\top}\widetilde{\mathbf  f}_i^{(r)}
               =-\,[\widetilde{\mathbf f}_i^{(k)}]^{\!\top}A_i\widetilde{\mathbf  f}_i^{(r)}$,
which is symmetric in $(k,r)$. Then
\begin{equation}
\begin{array}{cc}
\mathbf p_i^{\!\top}G_i\,\mathbf p_i
= -\sum_{k,r}p_{i,k}\,[\widetilde{\mathbf  f}_i^{(k)}]^{\!\top}A_i\widetilde{\mathbf  f}_i^{(r)}\,p_{i,r}
&= -\Bigl(\sum_{r}p_{i,r}\widetilde{\mathbf  f}_i^{(r)}\Bigr)^{\!\top}
    A_i
    \Bigl(\sum_{s}p_{i,s}\widetilde{\mathbf  f}_i^{(s)}\Bigr) \\ 
&    
    \overset{\eqref{eq:vel_pressure}}{=} -\,\mathbf U_i^{\!\top}A_i\,\mathbf U_i,
\end{array}
\end{equation}
which is the relation \eqref{eq:local_energy_identity}. 
Nonpositivity is immediate since $A_i$ is SPD.
\end{proof}

\begin{lemma}[Constant-traction kernel]\label{lem:const_kernel}
Let $\mathbf 1_i\in\mathbb R^{m_i}$ be the all-ones vector. Then
$G_i\mathbf 1_i=\mathbf 0$, and more generally,
$\mathbf p_i^{\!\top}G_i\,\mathbf p_i=0$ if and only if the interface traction is uniform,
i.e.\ $p_{i,1}=\cdots=p_{i,m_i}$.
\end{lemma}

\begin{proof}
Set $p_{i,r}\equiv p_{\mathrm{const}}$. The right-hand side of the local
Stokes problem \eqref{eq:stokes_local_matrix} is then
$\sum\limits_{r=1}^{m_i} p_{i,r}\mathbf f_i^{(r)}=p_{\mathrm{const}}\sum\limits_{r=1}^{m_i}\mathbf f_i^{(r)}$.
By the divergence theorem (noting that each velocity basis $\boldsymbol\phi_{i,j}$ vanishes on the solid wall $\partial \Omega_i\setminus (\cup_{r=1}^{m_i} \Gamma_{i,r})$) and the definition of $\mathbf f_i^{(r)}$,
\[
\sum_{r=1}^{m_i}\bigl[\mathbf f_i^{(r)}\bigr]_j
= -\sum_{r=1}^{m_i}\int_{\Gamma_{i,r}}\boldsymbol\phi_{i,j}\!\cdot \mathbf n_{i,r}\, dx
=
-\int_{\Omega_i}\nabla\!\cdot\boldsymbol\phi_{i,j}\,dx.
\]
Let $\pi_{i,\mathrm{const}}^h\in Q_i^h$ be the constant pressure function $\pi_{i,\mathrm{const}}^h\equiv p_{\mathrm{const}}$, with coefficient vector
$\mathbf P_{i,{\mathrm{const}}}$. By definition of $B_i$,
\[
\bigl[B_i^{\top}\mathbf P_{i,{\mathrm{const}}}\bigr]_j
= -\int_{\Omega_i}\pi^h_{i,\rm const}\,\nabla\!\cdot\boldsymbol\phi_{i,j}\,dx
= -p_{\rm const}\!\int_{\Omega_i}\nabla\!\cdot\boldsymbol\phi_{i,j}\,dx.
\]
Therefore, component-wise,
\[
\sum_{r=1}^{m_i} p_{\rm const}\,\mathbf f_i^{(r)} \;=\; B_i^{\top}\mathbf P_{i,{\mathrm{const}}}.
\]
Hence
$(\mathbf U_i,\mathbf P_i)=(\mathbf 0,\mathbf P_{i,{\mathrm{const}}})$ solves
\eqref{eq:stokes_local_matrix}, and by uniqueness (Lemma \ref{lem:local_wellposed}) $\mathbf U_i=\mathbf 0$.
Applying \eqref{eq:local_energy_identity} gives $\mathbf p_i^{\!\top}G_i\mathbf
p_i=0$.

For the converse implication ("$\mathbf p_i^{\!\top}G_i\mathbf p_i=0\Rightarrow$ uniform interface traction"), the proof  relies on two auxiliary results (Lemma \ref{lem:ker_char} and Lemma~\ref{lem:flux_surj}), and the full argument is deferred to \ref{app:aux}.

%if $\mathbf p_i^{\!\top}G_i\mathbf p_i=0$, then
%$\mathbf u_i^{h}=\sum_r p_{i,r}\widetilde f_i^{(r)}=\mathbf 0$
%(by \eqref{eq:local_energy_identity}). The first block row of
%\eqref{eq:stokes_local_matrix} with $\mathbf u_i^{h}=\mathbf 0$ gives
%$B_i^{\top}\pi_i^{h}=\sum_r p_{i,r}f_i^{(r)}$, which (by the definition of the
%face loads $f_i^{(r)}$) is precisely the statement that the imposed traction
%$-p_{i,r}\mathbf n$ equals $-\pi_i^{h}\mathbf n$ on every face. Testing with
%face-localized velocity functions shows $p_{i,r}=\pi_i^{h}|_{\Gamma_{i,r}}$
%for each $r$, hence all $p_{i,r}$ are equal.

Finally, set $p_{i,r}\equiv 1$. The corresponding discrete velocity is
$\displaystyle \mathbf U_i=\sum_{r=1}^{m_i}\widetilde{\mathbf f}_i^{(r)}$. By the above result of 
uniform interface implies zero velocity, $\mathbf U_i=\mathbf 0$.
For each row index $k$,
\[
\bigl[G_i\mathbf 1_i\bigr]_k
= \sum_{r=1}^{m_i} g_{i,kr}
= -\,[\mathbf f_i^{(k)}]^{\!\top}\!\Bigl(\sum_{r=1}^{m_i}\widetilde{\mathbf  f}_i^{(r)}\Bigr)
= -\,[\mathbf f_i^{(k)}]^{\!\top}\mathbf U_i
= 0.
\]
Hence $G_i\mathbf 1_i=\mathbf 0$.
\end{proof}

\begin{theorem}[Solvability and mass conservation of DD-PNM]
\label{thm:ddpnm}
With the local DtN operator $G_i$ defined in~\eqref{eq:dtn_local}
and the global Schur complement \\
$S=-\sum_{i=1}^{N}(R_i^{\rm u})^{\!\top}G_i^{\rm uu}R_i^{\rm u}$
introduced in \eqref{eq:global_schur},
the following properties of the DD-PNM hold: 
\begin{enumerate}\setlength{\itemsep}{4pt}
\item[\rm(i)] Positive definiteness.\;
      The global Schur complement
      $S$ is symmetric positive-definite (SPD).
\item[\rm(ii)] Global solvability.\;
      The linear system $S\,\mathbf p=\mathbf F$
      (with the forcing vector $\mathbf F$ defined in~\eqref{eq:global_schur})
      has a unique solution for the unknown interface tractions
      $\mathbf p \in\mathbb R^{m}$.
\item[\rm(iii)] Mass conservation.\;
      The reconstructed velocity field
      $\mathbf u^{\mathrm{DD}}=\bigcup_{i}\mathbf u_i^{h}$
      satisfies both the  local and global mass balances:
      \[
        \sum_{r=1}^{m_i}Q_{i,r}=0 \quad(\forall\,i),\qquad
             \int_{\partial\Omega}\mathbf u^{\mathrm{DD}}\!\cdot\mathbf n\,dx = 0.
      \]
\end{enumerate}
\end{theorem}

\begin{proof}%[Proof of Theorem~\ref{thm:ddpnm}]

\noindent
{(i) Symmetry and positive definiteness of $S$.} 
For any split $\mathbf p_i=(\mathbf p_i^{\mathrm u},\mathbf p_i^{\mathrm k})$ with
$\mathbf p_i^{\mathrm k}=\mathbf 0$, Lemma~\ref{lem:energy_identity} implies
\[
\mathbf p_i^{\!\top}G_i\,\mathbf p_i
=(\mathbf p_i^{\mathrm u})^{\!\top}G_i^{\rm uu}\mathbf p_i^{\mathrm u}
=-\,\mathbf U_i^{\top}A_i\,\mathbf U_i\le 0,
\]
so $G_i^{\rm uu}$ is negative semi-definite and contributes a
positive semi-definite block to
$S=-\sum_i(R_i^{\rm u})^{\!\top}G_i^{\rm uu}R_i^{\mathrm u}$.
Strict positivity holds unless $\mathbf U_i\equiv 0$, which by
Lemma \ref{lem:const_kernel} would require $p_{i,r}$ to be uniform on
\emph{all} interfaces. The only uniform choice compatible with
$\mathbf p_i^{\rm k}=0$ is $\mathbf p_i^{\rm u}=\mathbf 0$; hence $S$ is SPD.

\medskip
\noindent
{(ii) Global solvability.} Because $S$ is SPD, the linear system
$S\,\mathbf p=\mathbf F$ 
has exactly one solution.

\medskip
\noindent
{(iii) Local and global mass conservation.}
The discrete incompressibility constraint reads:
$\int_{\Omega_i} q_h\,\nabla\!\cdot\mathbf u_i^h\,dx=0$ for all
$q_h\in Q_i^h$. Choosing $q_h\equiv 1\in Q_i^h$ yields
$\int_{\Omega_i}\nabla\!\cdot\mathbf u_i^h\,dx=0$. By the
divergence theorem and the no-slip condition on solid walls,
\[
0=\int_{\Omega_i}\nabla\!\cdot\mathbf u_i^h\,dx
=\sum_{r=1}^{m_i}\int_{\Gamma_{i,r}}\mathbf u_i^h\!\cdot \mathbf n_{i,r}\,dx,
\]
which establishes the local mass balance.

For the global balance note that each internal interface contributes
equal and opposite flux to the neighbouring subdomains; thus all
internal fluxes cancel and only the flux through the
\emph{external} boundary $\partial\Omega$ remains,
yielding
$\displaystyle\int_{\partial\Omega}\mathbf u^{\mathrm{DD}}\!\cdot\mathbf n\,dx=0$.
\end{proof}

\section{Numerical results}\mbox{}

\label{sec:num_results}

\subsection{model problem}

%\textcolor{blue}{????I should describe the computational domain}

\label{sec:num_model_pb}

\begin{figure}[htbp!]
		\centering
%		\subfloat[$\vert \mathbf u \vert$]{\includegraphics[scale=0.34]{figures/FEM_case1_u.eps}}
%		\hfill
%		\subfloat[$p$]{\includegraphics[scale=0.34]{figures/FEM_case1_p.eps}}
\includegraphics[scale=0.34]{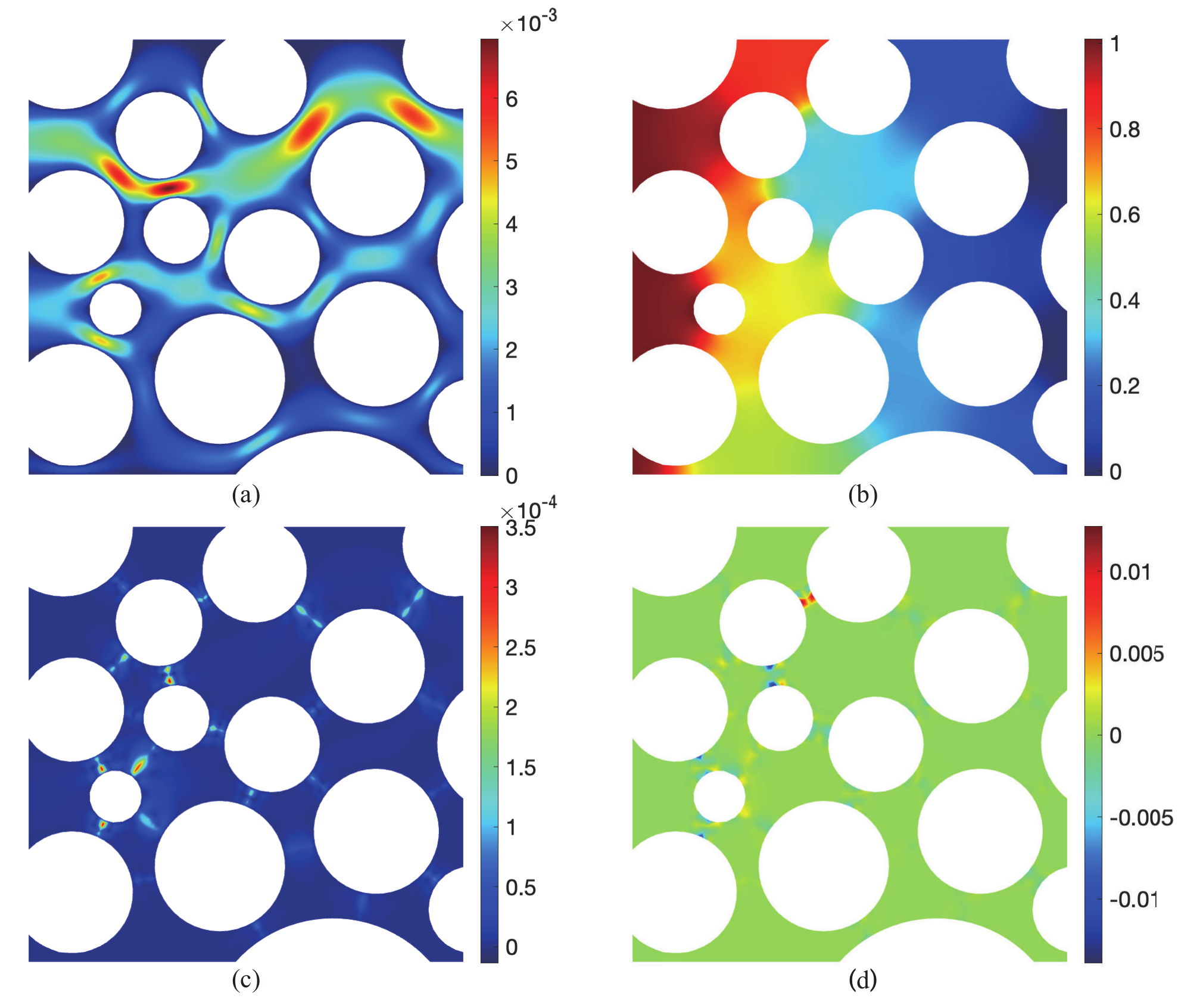}
		\caption{Model problem; (a)-(b) FE solution: (a) $\vert \mathbf u \vert$ and (b) $p$; (c)-(d) point-wise errors between DD-PNM and FE solutions: (c) $\vert \mathbf u \vert$ and (d) $p$.}
		\label{fig:FEM_sol_case1}
	\end{figure}

\begin{figure}[htbp!]
		\centering
%		\subfloat[$\varepsilon=0.1$]{\includegraphics[scale=0.34]{figures/DDPNM_mesh_shift0.1.eps}}
%		\hfill
%		\subfloat[$\varepsilon=0.2$]{\includegraphics[scale=0.34]{figures/DDPNM_mesh_shift0.2.eps}}
%		\\
%		\vspace*{-0.4cm}
%		\subfloat[$p$]{\includegraphics[scale=0.34]{figures/DDPNM_shift0.1_err_p.eps}}
%		\hfill
%		\subfloat[$p$]{\includegraphics[scale=0.34]{figures/DDPNM_shift0.2_err_p.eps}}
\includegraphics[scale=0.22]{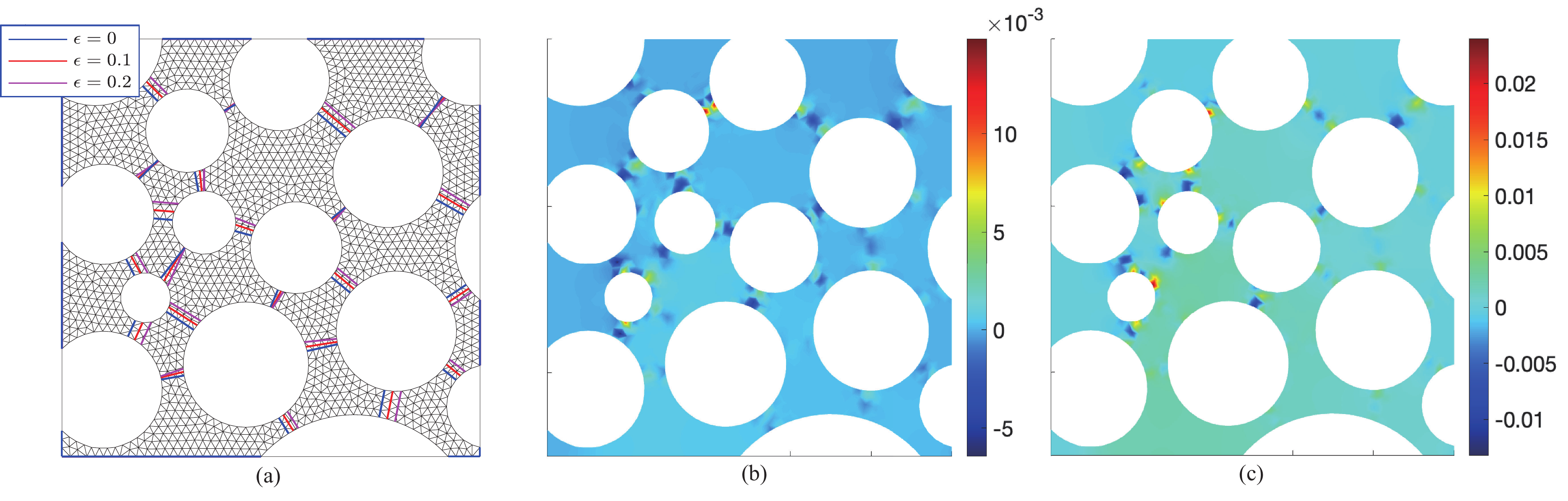}
		\caption{Model problem; DD-PNM solution with perturbed interfaces: (a) interface positions for perturbation magnitudes $\varepsilon=0.1$ and $\varepsilon=0.2$; (b)-(c) corresponding pressure error fields relative to the FE reference solution. }
				\label{fig:insensibility_case1}
	\end{figure}

We revisit the two-dimensional Stokes flow configuration introduced in Section~\ref{sec:model_problem}.
The geometry consists of a rectangular pore domain, $[0,5]\times[0,5]$, with randomly packed circular solid inclusions. The top and bottom boundaries are treated as no-slip walls, 
and flow is driven form left to right by imposing a pressure drop from inlet ($p_{\rm in}=1$) to the outlet ($p_{\rm out}=0$).
To assess the performance of the proposed DD-PNM method, we use a  finite-element (FE) solution as a reference. 
Figure~\ref{fig:FEM_sol_case1} presents the FE solution and the point-wise errors between the DD-PNM and the FE reference. Panels (a) and (b) illustrate the FE velocity magnitude and pressure fields, respectively, while panels (c) and (d) show the corresponding point-wise error fields. 
%These errors are concentrated within narrow regions near the throat interfaces, reflecting the modeling simplifications introduced by the DD-PNM interface conditions.
%Figure \ref{fig:FEM_sol_case1} shows the velocity and pressure fields obtained from the FE simulation. Panel (a) depicts the velocity magnitude, and panel (b) presents  the pressure field.
%Figure~\ref{fig:DDPNM_FEM_err_case1}
%presents the point-wise errors in velocity and pressure between the DD-PNM and FE solutions.
  Both errors are confined to narrow regions near the throat interfaces. %, confirming the high accuracy of the proposed method.
For the FE solution using the coarse mesh (see Figure \ref{fig:insensibility_case1}(a)), the resulting linear system has $12765$ degrees of freedom (DOFs). In contrast, the DD-PNM approach solves subdomain problems in parallel, each ranging from $282$ to $1486$ DOFs, coupled by a small Schur complement system of size $22\times 22$. Thus, DD-PNM achieves significant computational efficiency.

%	\begin{figure}[htbp!]
%		\centering
%		\subfloat[$\vert \mathbf u \vert$]{\includegraphics[scale=0.34]{figures/DDPNM_FEM_err_u_case1.eps}}
%		\hfill
%		\subfloat[$p$]{\includegraphics[scale=0.34]{figures/DDPNM_FEM_err_p_case1.eps}}
%		\caption{Model problem; DD-PNM versus FE: point-wise error fields.}
%		\label{fig:DDPNM_FEM_err_case1}
%	\end{figure}

% DD-PNM vs PNM, q \\

To assess robustness with respect to interface location, we introduce perturbations in the placement of each interface. As described in Section \ref{sec:ddpnm}, interfaces in the baseline model are defined as lines connecting neighboring solid centers. Here, we perturb the endpoints of these lines: rather than using the original endpoint $\mathbf x_0$, we select perturbed point $\mathbf x_0 + \mathbf y$, where $\mathbf y$ is a random vector uniformly sampled from the ball $B(0,r\varepsilon)=\{\mathbf{x}:\|\mathbf x\|\leq r\varepsilon\}$, with $r$ being a random number in the interval $[0,1]$. 
Figure \ref{fig:insensibility_case1}(a) illustrates the perturbed interfaces for perturbation magnitudes $\varepsilon=0.1$ and $\varepsilon=0.2$. For simplicity, only the original (unperturbed) mesh is shown; in practice, a new mesh would be generated according to each perturbed interface position. Panels (b) and (c) depict the corresponding pressure errors relative to the FE reference solution. 
In both cases, the errors remain small and spatially localized, indicating that the DD-PNM method is robust with respect to variations in interface placement.

\begin{figure}[htbp!]
		\centering
%		\subfloat[$h_0=0.1$]{\includegraphics[scale=0.23]{figures/FEMmesh_case1_local.eps}}
%		\hfill
%		\subfloat[$h_0=0.05$]{\includegraphics[scale=0.23]{figures/FEMmesh_case1_refine1.eps}}
%		\hfill
%		\subfloat[$h_0=0.025$]{\includegraphics[scale=0.23]{figures/FEMmesh_case1_refine2.eps}}
%		\hfill
%		\vspace*{-0.3cm}
%		\\
%		\subfloat[$\vert \mathbf u \vert$]{\includegraphics[scale=0.34]{figures/FEM_case1_refine1_u.eps}}
%		\hfill
%		\subfloat[$p$]{\includegraphics[scale=0.34]{figures/FEM_case1_refine1_p.eps}}
%		\\
%		\vspace*{-0.4cm}
%		\subfloat[$\vert \mathbf u \vert$]{\includegraphics[scale=0.34]{figures/FEM_case1_refine2_u.eps}}
%		\hfill
%		\subfloat[$p$]{\includegraphics[scale=0.34]{figures/FEM_case1_refine2_p.eps}}
		\includegraphics[scale=0.22]{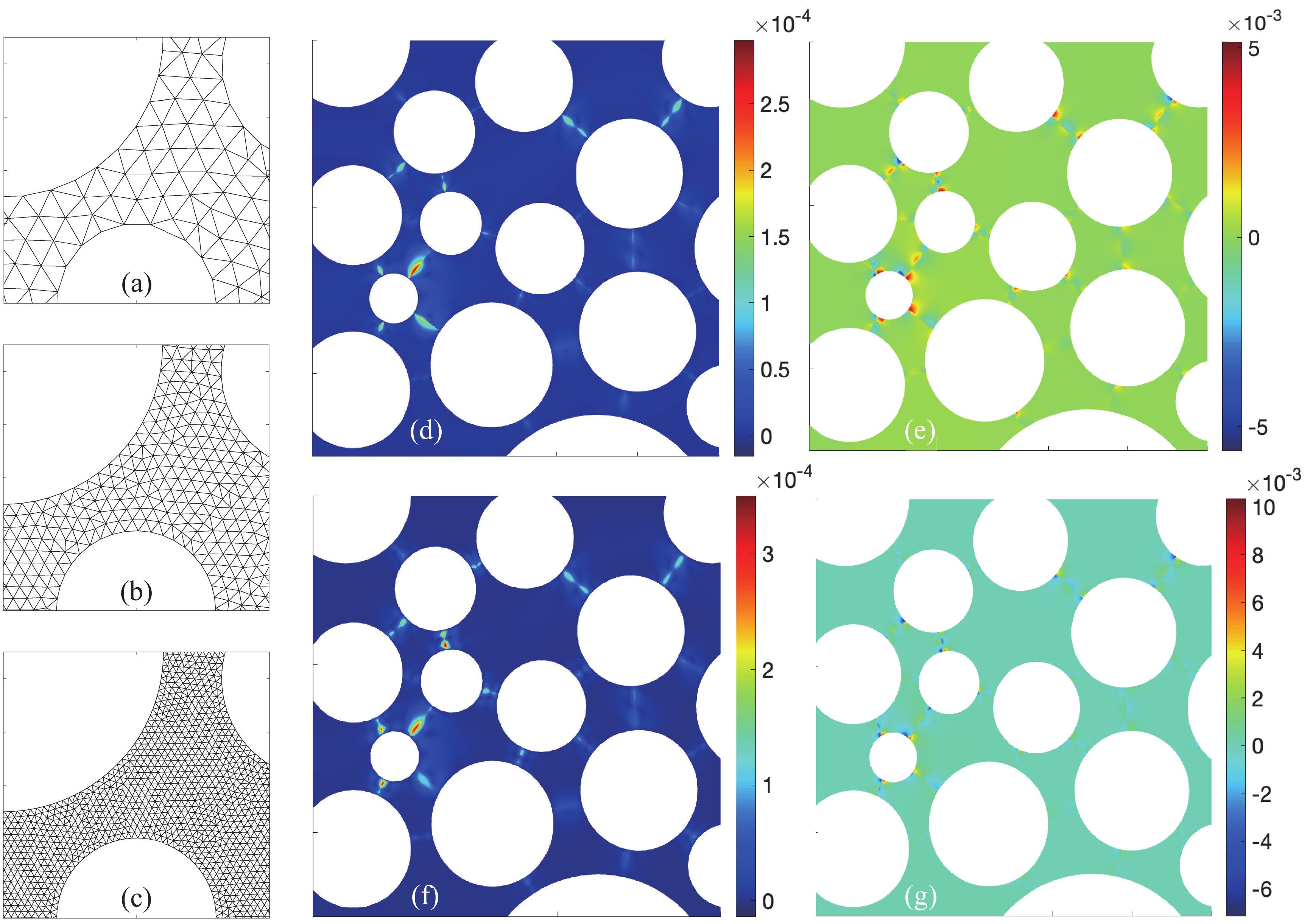}
		\caption{Model problem; DD-PNM versus FE solutions on refined meshes. Panels (a)-(c) show zoomed-in views of the original mesh ($h_0=0.1$) and two refined meshes ($h_0=0.05$ and $h_0=0.025$). Panels (d)-(e) present velocity and pressure error fields for $h_0=0.05$; (f)-(g) velocity and pressure error fields for $h_0=0.025$.}
		\label{fig:refined_mesh_case1}
	\end{figure}

We next evaluate the accuracy of DD-PNM on two refined FE meshes with characteristic sizes $
h_0=0.05$ and 
$ h_0=0.025$, respectively. 
For each refinement, we compare the resulting DD-PNM velocity and pressure fields against the corresponding FE solutions.
Figure \ref{fig:refined_mesh_case1}
displays zoomed-in views of the original and refined meshes for clearer visualization (panels (a)–(c)), alongside the corresponding velocity and pressure error fields.
It is observed that the point-wise error of the velocity magnitude and the point-wise error of the pressure varies  only slightly when the mesh is refined (from $h_0=0.1$ to $h_0 = 0.05$ and $h_0 = 0.025$). 
This observation aligns with the error estimate discussed in Appendix \ref{sec:error_estimate} (see Remark \ref{remark:error_estimate}): the error does not decrease unless the interface approximation space is enriched. 

%Efficient numerical methods with improve accuracy will be future work. 
%The DD-PNM solution converges rapidly toward the FE reference as the local pore meshes are refined.

	Table \ref{table:field_error_case1} summarizes the relative $H^1$-velocity and $L^2$-pressure errors for the following five scenarios:  run	1 (origin case $h_0 = 0.1$), 
run	2-3 (refined meshes, $h_0 = 0.05,\, 0.025$);
run	4-5 (interface locations perturbed with $\varepsilon = 0.1,\, 0.2$). 
	In all cases, the relative 
	 error remains at  $\mathcal{O}(10^{-3})$, demonstrating the accuracy and robustness of the DD-PNM.
	\begin{table}[htbp!]
		\centering
		\caption{Model problem; relative $H^1\times L^2$ errors for the velocity-pressure pair.}
			\begin{tabular}{|c|c|c|c|c|c|}
			\hline
				run & 1 & 2 & 3 & 4 & 5  \\
				\hline
				rel error & 0.0017 & 0.0013 & 0.0017 & 0.0027 & 0.0059	\\
				\hline
			\end{tabular}
		\label{table:field_error_case1}
	\end{table}

\begin{figure}[htbp!]
		\centering
		\includegraphics[scale = 0.3]{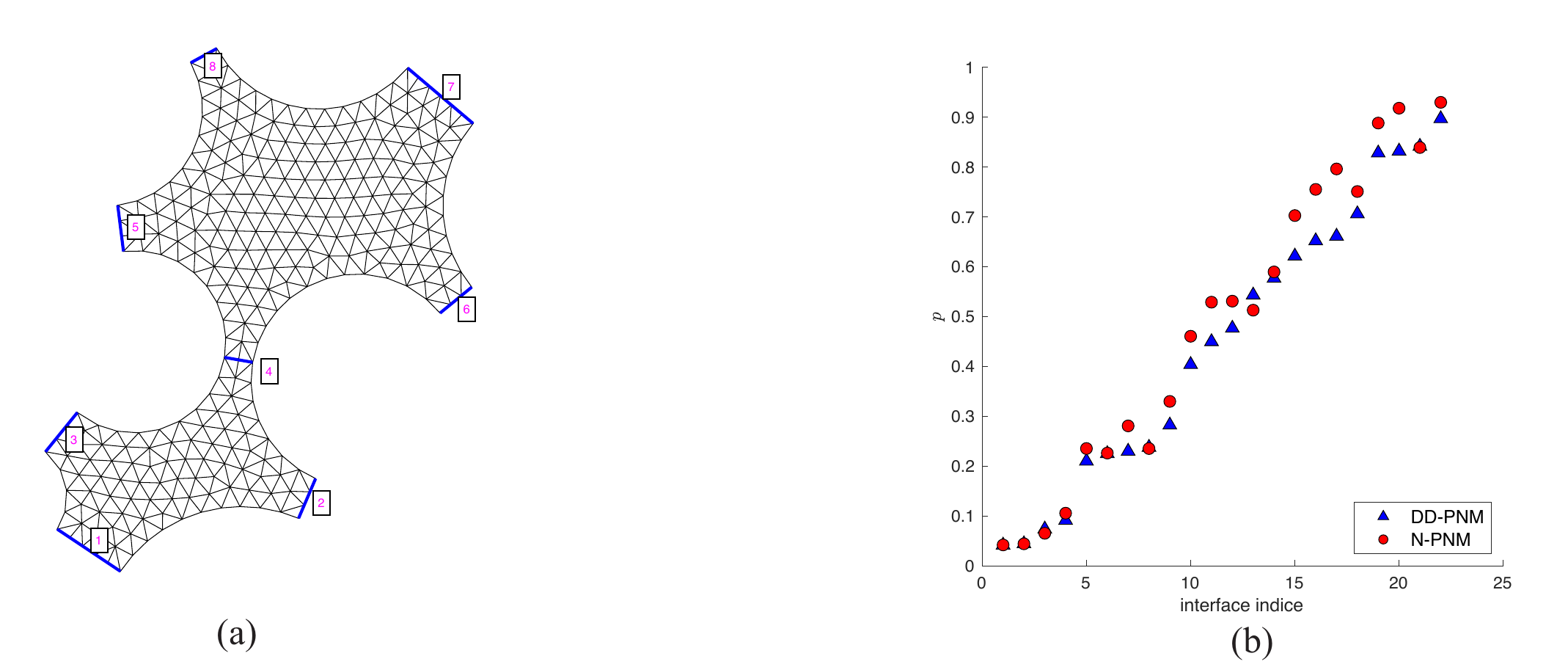}
		\caption{Model problem; (a) two adjacent pore regions used for the local mass-balance test; (b) comparison of interface tractions computed by DD-PNM and N-PNM. }
		\label{fig:loc_region_mesh}
	\end{figure}

We now verify that the DD-PNM preserves local mass balance across pore regions.
Figure \ref{fig:loc_region_mesh}(a)  shows two adjacent subdomains and their shared interfaces.
Table \ref{table:flux} reports the computed volumetric flux across each labeled interface. It can be directly verified that, 
for each region, the net flux sums to zero, and the fluxes across the shared interface cancel out. This confirms that the proposed method strictly enforces discrete mass  conservation.

	\begin{table}[htbp!]
		\centering
		\scriptsize % Smaller font size
		\setlength{\tabcolsep}{3pt} % Reduce column spacing
		\caption{Model problem; volumetric fluxes across labeled interfaces.}
		\resizebox{1\textwidth}{!}{
	\begin{tabular}{|c|c|c|c|c|}
		\hline
		interface & 1 & 2 & 3 & 4  \\
		\hline
		Q & -3.4420 & 6.5636 & -5.5008 & 2.3792 	\\
		\hline
	\end{tabular}
	\hfill
		\begin{tabular}{|c|c|c|c|c|c|}
		\hline
		interface & 4 & 5 & 6 & 7 & 8  \\
		\hline
		Q & -2.3792 & -12.4551 & 2.9572 & 16.2517 & -4.3746	\\
		\hline
	\end{tabular}
}
\label{table:flux}
	\end{table}

We denote the C-PNM recovered from DD-PNM, as introduced in Section \ref{sec:recovering_CPNM}, as the New-PNM (N-PNM). Table \ref{table:CpnmvsNpnm} compares
 the pore-body pressures obtained using N-PNM and C-PNM, showing reasonable agreement between the two methods.
%The two methods show reasonable agreement, with noticeable differences reflecting the distinct approaches used in determining half-throat conductivities.
Furthermore, using the half-throat conductivities derived in N-PNM, we  compare the interface normal tractions computed by N-PNM and DD-PNM. 
Figure \ref{fig:loc_region_mesh}(b) illustrates these interface tractions from both methods. For clarity, interface indices have been sorted according to increasing DD-PNM traction values. The close agreement between N-PNM and DD-PNM tractions demonstrates the effectiveness of our proposed approach to recovering C-PNM from DD-PNM.

		\begin{table}[htbp!]
		\centering
		\scriptsize % Smaller font size
		\setlength{\tabcolsep}{3pt} % Reduce column spacing
		\caption{Model problem; comparison of pore-body pressures computed by C-PNM and N-PNM.}
		\resizebox{1\textwidth}{!}{
			\begin{tabular}{|c|c|c|c|c|c|c|c|c|c|c|c|c|c|c|c|c|c|}
				\hline
				pore body & $1$ & $2$ & 3 & 4 & 5 & 6 & 7 & 8 & 9 & 10 & 11 & 12 & 13 & 14 & 15 & 16 & 17  \\
				\hline
				$p_\text{N-PNM}$  & 0.5922 & 0.3487 & 0.1990 & 1 & 0.7480 & 0.7033 & 0.3594 & 0 & 0.8448 & 0.7091 & 0.3931 & 0.1196 & 0.1134 & 0 & 1 & 0.7907 & 0.1248 	\\
				\hline
				$p_\text{C-PNM}$  &0.6156 & 0.1942 & 0.1714 & 1 & 0.7130 & 0.7070 & 0.1870 & 0 & 0.7139 & 0.8608 & 0.1344 & 0.0553 & 0.0533 & 0 & 1 & 0.9267 & 0.0547\\
				\hline
			\end{tabular}
		}
		\label{table:CpnmvsNpnm}
	\end{table}

%	\begin{figure}[htbp!]
%		\centering
%		\includegraphics[scale = 0.3]{figures/scatter_case1.eps}
%		\caption{Model problem; comparison of interface tractions computed by DD-PNM and N-PNM. }
%		\label{fig:scatter_case1}
%	\end{figure}

\subsection{Realistic geometry}
\label{sec:num_real_geo}

\begin{figure}[htbp!]
		\centering
		\includegraphics[width=13cm]{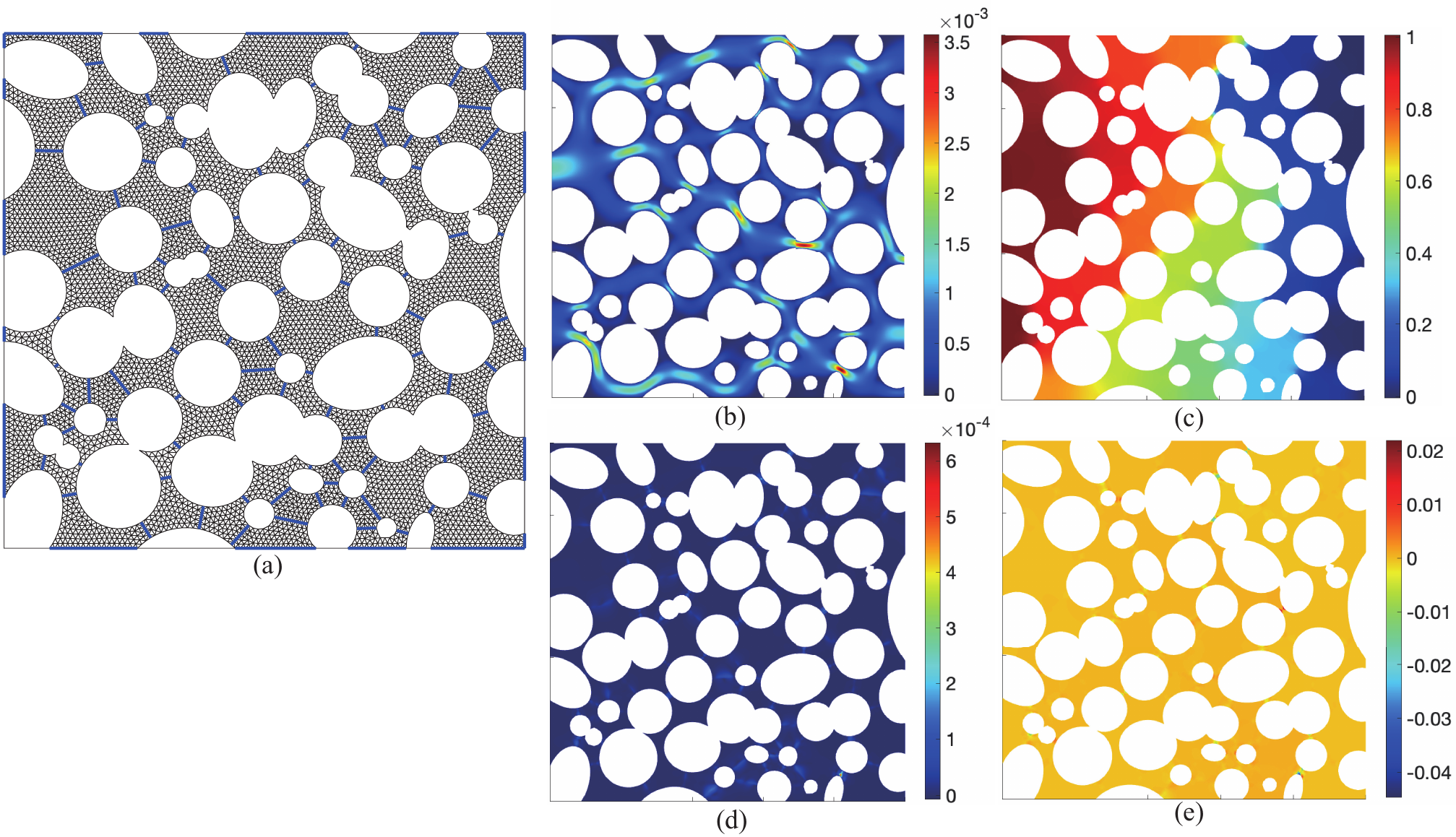}
		\caption{Realistic geometry; (a) computational geometry and mesh; (b)-(c) FE solution: (b) $\vert \mathbf u \vert$ and (c) $p$; (d)-(e) point-wise erros between DD-PNM and FE solutions: (d) $\vert \mathbf u \vert$ and (e) $p$. }
		\label{fig:DDPNMmesh_case2}
	\end{figure}

We now consider a more realistic geometry, depicted in Figure \ref{fig:DDPNMmesh_case2}(a), in which the solid particles include elliptical shapes and meta-balls within a rectangular domain $[0,10]\times[0,10]$. The figure also illustrates the corresponding domain decomposition and the generated computational meshes. The boundary conditions applied are identical to those used in the previous test case described in Section \ref{sec:num_model_pb}. 
For this geometry, we first apply maximal ball algorithm to identify the approximate location of interfaces; we then precisely define the interface positions using the minimal distance between adjacent particles. 

Figure \ref{fig:DDPNMmesh_case2}(b)-(c) presents the finite-element solutions for the velocity magnitude and pressure fields associated with this geometry. 	Panels (d)-(e) shows the corresponding point-wise velocity and pressure errors between the DD-PNM and FE solutions. 
	Similar to the previous test case, errors are primarily localized near the throat interfaces. 
	The relative error measured in $H^1\times L^2$  norm for the velocity-pressure pair is $0.0013$, demonstrating the high accuracy of the proposed method. We report the computational efficiency of the DD-PNM approach. While the FE method solves a system with $48944$ DOFs, the DD-PNM method solves smaller local problems in parallel ranging from $163$ to $3201$ DOFs, coupled with a Schur complement system of $91\times 91$.

%		\begin{figure}[htbp!]
%		\centering
%		\subfloat[$\vert \mathbf u \vert$]{\includegraphics[scale=0.34]{figures/FEM_case2_u.eps}}
%		\hfill
%		\subfloat[$p$]{\includegraphics[scale=0.34]{figures/FEM_case2_p.eps}}
%		\caption{Realistic geometry; FE solution.}
%		\label{fig:FEM_sol_case2}
%	\end{figure}

%	\begin{figure}[htbp!]
%		\centering
%		\subfloat[$\vert \mathbf u \vert$]{\includegraphics[scale=0.34]{figures/DDPNM_FEM_err_u_case2.eps}}
%		\hfill
%		\subfloat[$p$]{\includegraphics[scale=0.34]{figures/DDPNM_FEM_err_p_case2.eps}}
%		\caption{Realistic geometry; DD-PNM versus FE: point-wise error fields.}
%		\label{fig:DDPNM_FEM_err_case2}
%	\end{figure}

Finally, we compare interface tractions obtained using N-PNM (see Section \ref{sec:num_model_pb}) and DD-PNM. The results align closely, further validating our approach to recovering C-PNM from DD-PNM.

	\begin{figure}[htbp!]
		\centering
		\includegraphics[scale = 0.3]{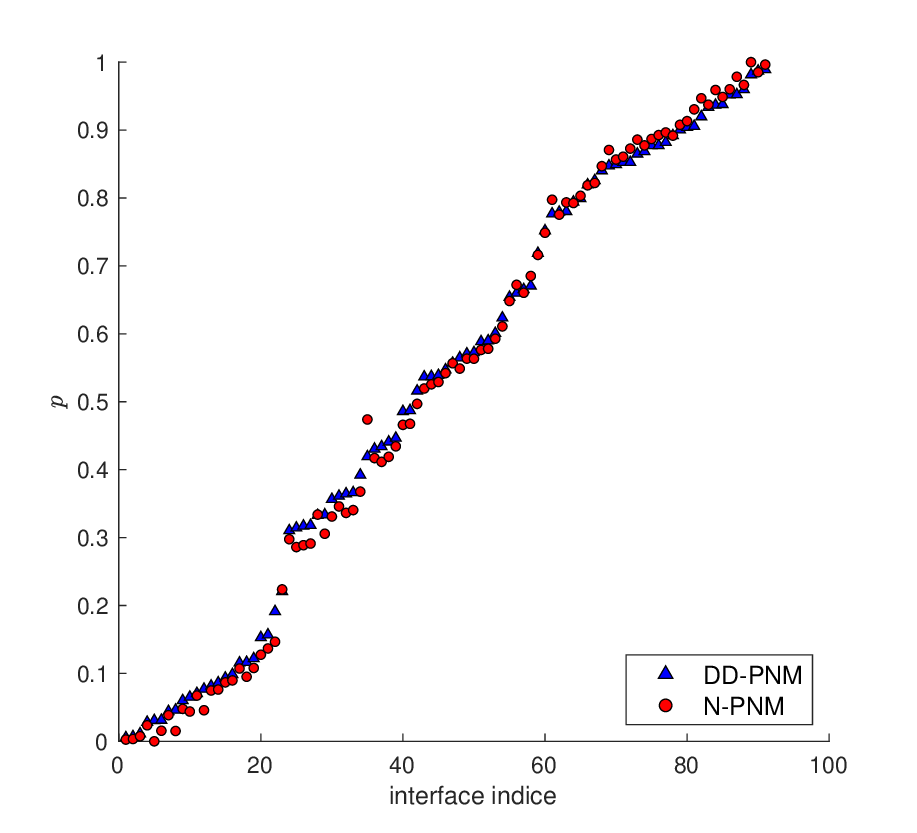}
		\caption{Realistic geometry; comparison of interface tractions between DD-PNM and N-PNM.}
		\label{fig:}
	\end{figure}

\section{Conclusion}

We introduced and analyzed a domain-decomposition \\
pore-network method (DD-PNM) that combines finite-element accuracy inside pore subdomains while retaining computational efficiency through sparse global interface coupling. The method rigorously enforces mass conservation and guarantees solvability of the resulting global system. 
By explicitly fitting half-throat conductivities to the computed Dirichlet-to-Neumann maps, we provided a rigorous approach to extract pore-throat conductivities as used in the classical network-based (C-PNM) pore-scale modeling framework; it can also be viewed as a constructive pathway connecting mesh-based (DD-PNM) and the C-PNM pore-scale modeling frameworks.
%By explicitly fitting half-throat conductivities to computed Dirichlet-to-Neumann maps, we provided a constructive pathway connecting mesh-based (DD-PNM) and classical network-based (C-PNM) pore-scale modeling frameworks.
Numerical results confirmed the accuracy and effectiveness of the proposed DD-PNM approach. Future work will focus on extending interface modeling strategies, such as Robin conditions, integrating non-Newtonian and multiphase flow physics, and evaluating performance on large-scale realistic micro-CT image datasets.

%\section*{Acknowledgements}

%Although cheap and robust, C-PNM accuracy degrades in complex throats where the circular-tube assumption is violated.  
%This motivates the hybrid approach proposed in Section~\ref{sec:ddpnm}.  
%Although C-PNM is remarkably efficient, its simplifying assumptions
%motivate us to seek a high-fidelity reference
%against which improvements can be quantified.

%To obtain benchmark solutions we solve the full Stokes problem on $\Omega$:

%The symmetric positive-definite matrix $G_i$ acts as a
%high-fidelity \emph{transmissibility tensor} that generalises the
%Hagen–Poiseuille coefficient. ????

%We note that for this simplified configuration, the PNM method is simple, for complex geometry, we may require complicated algorithm. Also we show our method is superior to the PNM method, advanced PNM method may improve the result, but we emphasize that our result is comparable to FEM results.  ????

\appendix

\section{Auxiliary results for Section~\ref{sec:properties_DD-PNM}}\label{app:aux}

\begin{lemma}[Kernel characterization]\label{lem:ker_char}
Let $F_i=[\,\mathbf f_i^{(1)}\ \cdots\ \mathbf f_i^{(m_i)}\,]$ be the matrix of the Neumann load vectors and let $G_i$ be the local DtN matrix defined in Section \ref{sec:ddpnm}. 
Then
\[
\mathbf{p}_i\in\ker(G_i)\iff F_i\mathbf{p}_i\in\mathrm{Range}(B_i^\top). 
\]
\end{lemma}

\begin{proof}

Fix $\mathbf p_i=[p_{i,1},\dots,p_{i,m_i}]^\top$. The local system \eqref{eq:stokes_local_matrix} can be written as
\begin{equation}
\label{eq:SPS}
\begin{bmatrix}A_i & B_i^\top\\[2pt] B_i & 0\end{bmatrix}
\begin{bmatrix}\mathbf U_i\\[2pt]\mathbf P_i\end{bmatrix}
=
\begin{bmatrix}F_i\,\mathbf p_i\\[2pt]0\end{bmatrix}.
\end{equation}
($\Longrightarrow$) If $\mathbf p_i\in\ker(G_i)$, then by the energy identity (Lemma \ref{lem:energy_identity}) we have
$\mathbf U_i=\mathbf 0$. Returning to the first block row of \eqref{eq:SPS} gives
\[
B_i^\top \mathbf P_i = F_i\,\mathbf p_i
\quad\Longrightarrow\quad
F_i\,\mathbf p_i  \in \mathrm{Range}(B_i^\top).
\]
($\Longleftarrow$) If $F_i\,\mathbf p_i\in\mathrm{Range}(B_i^\top)$, write $F_i\,\mathbf p_i=B_i^\top \widehat{\mathbf P}$ for some $\widehat{\mathbf P}$.
Then $(\widehat{\mathbf U},\widehat{\mathbf P})=(\mathbf 0,\widehat{\mathbf P})$ satisfies \eqref{eq:SPS}. By uniqueness,
$\mathbf U_i=\mathbf 0$, and Lemma \eqref{lem:energy_identity}
gives 
$\mathbf p_i^\top G_i\,\mathbf p_i =  0$.
Since $G_i$ is symmetric negative semidefinite
(Lemma~\ref{lem:energy_identity}), $\mathbf p_i^\top G_i\,\mathbf p_i=0$ is equivalent to $G_i\mathbf p_i=\mathbf 0$, i.e., $\mathbf p_i\in\ker(G_i)$.
\end{proof}

\begin{lemma}[Flux-surjectivity onto zero-sum vectors]\label{lem:flux_surj}
Let $\mathcal{Z}_i:=\ker(B_i)$ and define $\Phi:\mathcal{Z}_i\to\mathbb R^{m_i}$ by
\begin{equation}
\Phi(\mathbf V):=\Bigl(\ \int_{\Gamma_{i,1}}\mathbf v_h\!\cdot\!\mathbf n_{i,1}\,dx,\ \dots,\ \int_{\Gamma_{i,m_i}}\mathbf v_h\!\cdot\!\mathbf n_{i,m_i}\,dx\ \Bigr)^{\!\top},
\label{eq:flux_map}
\end{equation}
where $\mathbf{v}_{h}$ is the FE function with coefficient vector $\mathbf{V}.$
Then $\mathrm{Range}(\Phi)=\mathbb R^{m_i}_0:=\{\boldsymbol\alpha\in\mathbb R^{m_i}:\sum\limits_{r=1}^{m_i}\alpha_r=0\}$.
\end{lemma}

\begin{proof}
For $\mathbf V\in\ker(B_i)$ we have, by definition, 
$\int_{\Omega_i}q_h\nabla\cdot\mathbf{v}_h\:dx=0$ for all $q_h\in Q_i^h.$
Choosing $q_h\equiv 1\in Q_i^h$ gives
$\int_{\Omega_i}\nabla\cdot\mathbf{v}_h\:dx=0$. 
Using the divergence theorem and the homogeneous Dirichlet condition on  solid-wall part of $\partial\Omega_i:$
$$
\sum_{r=1}^{m_i}[\Phi(\mathbf{V})]_r\:=\:
\sum_{r=1}^{m_i}\int_{\Gamma_{i,r}}\mathbf{v}_h\cdot\mathbf{n}_{i,r}\:dx\:=\:\int_{\partial\Omega_i}\mathbf{v}_h\cdot\mathbf{n}\:dx\:=\:\int_{\Omega_i}\nabla\cdot\mathbf{v}_h\:dx\:=\:0.$$
Hence
$$\mathrm{Range}(\Phi)\subseteq\mathbb{R}_0^{m_i}:=\{\boldsymbol{\alpha}\in\mathbb{R}^{m_i}:\sum_r\alpha_r=0\}.$$
Equip $\mathcal{Z}_i$ with $(\mathbf{U},\mathbf{V})_{A_i}:=\mathbf{U}^\top A_i\mathbf{V}$ (an inner product since $A_i$ is SPD). 
For any $\boldsymbol{\alpha}\in\mathbb{R}_0^{m_i}$, define the linear functional on $\mathcal{Z}_i:$
$$\mathcal{L}_{\boldsymbol{\alpha}}(\mathbf{V}):=\boldsymbol{\alpha}^\top\Phi(\mathbf{V}).$$
By the Riesz representation theorem, there exists a unique $\mathbf{W}\in \mathcal{Z}_i$ such that
$$\mathbf{W}^\top A_i\mathbf{V}=\mathcal{L}_{\boldsymbol{\alpha}}(\mathbf{V})\quad\forall\:\mathbf{V}\in \mathcal{Z}_i.$$
Evaluating at $\mathbf{V}=\mathbf{W}$ yields
$$
\boldsymbol{\alpha}^{\top}\Phi(\mathbf{W})=\mathcal{L}_{\boldsymbol{\alpha}}(\mathbf{W}) = 
\mathbf{W}^{\top}A_{i}\mathbf{W}>0\quad\mathrm{if~}\boldsymbol{\alpha}\neq\boldsymbol{0}.$$
Thus the linear operator $M:\mathbb{R}_0^{m_i}\to\mathbb{R}_0^{m_i}$ defined by
$$M(\boldsymbol{\alpha}):=\Phi(\mathbf{W})\quad(\mathrm{with~}\mathbf{W}\text{ as above})$$
is SPD on $\mathbb{R}_0^{m_i}$ (since $\boldsymbol{\alpha}^\top M\boldsymbol{\alpha}=\mathbf{W}^\top A_i\mathbf{W}>0$ for $\boldsymbol{\alpha}\neq0$). Therefore $M$ is
invertible on $\mathbb{R}_0^{m_i}$, and in particular it is onto $\mathbb{R}_0^{m_i}$. That is, Range$( \Phi ) = \mathbb{R} _0^{m_i}.$
\end{proof}

\begin{proof}[Proof of the converse implication in Lemma~\ref{lem:const_kernel}]

Assume $\mathbf p_i^{\!\top}G_i\mathbf p_i=0$. \\
By Lemma \ref{lem:ker_char}, this implies $F_i\mathbf p_i\in\mathrm{Range}(B_i^\top)$; 
hence there exists $\widehat{\mathbf P}$ such that $F_i\mathbf p_i=B_i^\top\widehat{\mathbf P}$.
Let $\mathbf V\in \mathcal{Z}_i=\ker(B_i)$ and let $\mathbf v_h$ be the FE velocity corresponding to $\mathbf V$. Taking the Euclidean inner product with $\mathbf V$ gives
$$
\mathbf V^\top F_i\mathbf p_i
=\mathbf V^\top B_i^\top\widehat{\mathbf P}
=(B_i\mathbf V)^\top\widehat{\mathbf P}
=0.
$$
By the definition of the Neumann load vectors,
$$
\sum_{r=1}^{m_i} p_{i,r}\,\mathbf V^\top\mathbf f_i^{(r)}
=-\sum_{r=1}^{m_i}p_{i,r}\int_{\Gamma_{i,r}}\mathbf v_h\!\cdot\!\mathbf n_{i,r}\,dx
=0
\qquad\forall\,\mathbf V\in \mathcal{Z}_i.
$$
Equivalently, $\mathbf p_i$ annihilates the range of the  map
$\Phi$ (see \eqref{eq:flux_map}). 
By Lemma \ref{lem:flux_surj}, $\mathrm{Range}(\Phi)=\mathbb R^{m_i}_0:=\{\boldsymbol\alpha\in\mathbb R^{m_i}:\sum_r \alpha_r=0\}$.
Therefore $\mathbf p_i$ is orthogonal to $\mathbb R^{m_i}_0$, i.e.,
$\mathbf p_i\in(\mathbb R^{m_i}_0)^\perp=\mathrm{span}\{\mathbf 1_i\}$.
Hence the interface traction is uniform, which completes the converse implication.
\end{proof}

\section{Error estimate of the DD-PNM}
\label{sec:error_estimate}

We derive error estimate for the DD-PNM approximation with respect to the
classical nonoverlapping domain decomposition (C-DD) formulation of the discrete Stokes problem. We first state the C-DD system. We then obtain the DD-PNM system as a Galerkin restriction of the C-DD formulation. This restricted form of the DD-PNM leads directly to the projection-type error estimate.

\subsection{Classical nonoverlapping DD}\label{sec:C-DD}

We follow the notation introduced in Section \ref{sec:ddpnm}. On internal interface $\Gamma_{i,r}$ of subdomain $\Omega_i$, the tractions decomposes into normal and tangential parts as
\[
\boldsymbol\sigma\mathbf n_{\vert{\Gamma_{i,r}}}
=\underbrace{\big(\mathbf n_{i,r}\cdot\!(\boldsymbol\sigma\mathbf n)_{\vert{\Gamma_{i,r}}}\big)}_{\lambda^n_{i,r}}\,\mathbf n_{i,r}
+\underbrace{\big(\mathbf t_{i,r}\cdot\!(\boldsymbol\sigma\mathbf n)_{\vert{\Gamma_{i,r}}}\big)}_{\lambda^t_{i,r}}\,\mathbf t_{i,r},
\]
with unit normal $\mathbf n_{i,r}$ and a unit tangent $\mathbf t_{i,r}$. 
We recall that in DD-PNM we enforce $\lambda^t_{i,r}\equiv 0$ and $\lambda^n_{i,r}\equiv p_{i,r}\in\mathbb R$ (one scalar per interface); in the  C-DD formulation, however, the traction is defined per mesh node on $\Gamma_{i,r}$: 
\begin{equation}
\begin{array}{l}
\lambda^n_{i,r} = \sum\limits_{\ell=1}^{m_{i,r}} (\boldsymbol\lambda^n_{i,r})_\ell\widehat\varphi_{i,r,\ell},
\quad
\lambda^t_{i,r} = \sum\limits_{\ell=1}^{m_{i,r}} (\boldsymbol\lambda^t_{i,r})_\ell\widehat\varphi_{i,r,\ell},
\end{array}
\end{equation}
where the vectors $\boldsymbol\lambda^n_{i,r}$ and $\boldsymbol\lambda^n_{i,r}$ collect normal/tangential traction degrees of freedom (DOFs) on $\Gamma_{i,r}$, 
$m_{i,r}$ is the number of interface DOFs on $\Gamma_{i,r}$, 
$\{\widehat\varphi_{i,r,\ell}\}$ is the finite element basis restricted to $\Gamma_{i,r}$. 

Let $n_{i,\rm u}$ (resp.\ $n_{i,\rm p}$) be the velocity (resp.\ pressure) DOFs on $\Omega_i$,
and $\{\boldsymbol\phi_j\}_{j=1}^{n_{i,\rm u}}$ the velocity basis. 
Define the
normal/tangential Neumann coupling blocks on each interface $\Gamma_{i,r}$ by
\begin{equation*}
\begin{array}{l}
\big(\widehat N^{\,n}_{i,r}\big)_{\ell j}
=\int_{\Gamma_{i,r}}\!\widehat\varphi_{i,r,\ell}\,(\boldsymbol\phi_j\!\cdot\mathbf n_{i,r})\,dx
\ \in\ \mathbb R^{m_{i,r}\times n_{i,\rm u}},\\
\big(\widehat N^{\,t}_{i,r}\big)_{\ell j}
=\int_{\Gamma_{i,r}}\!\widehat\varphi_{i,r,\ell}\,(\boldsymbol\phi_j\!\cdot\mathbf t_{i,r})\,dx
\ \in\ \mathbb R^{m_{i,r}\times n_{i,\rm u}}.
\end{array}
\end{equation*}
Then, the discrete Stokes system on each $\Omega_i$  reads
% ($\mathcal F_i^{\mathrm{int}}$ (resp.\ $\mathcal F_i^{\mathrm{bc}}$) denotes the sets of internal (resp.\ boundary) faces of $\Omega_i$)
\begin{equation}
\label{eq:local-stokes_dd}
\underbrace{\begin{bmatrix}A_i & B_i^\top\\ B_i & 0\end{bmatrix}}_{K_i}
\begin{bmatrix}\mathbf U_i\\ \mathbf P_i\end{bmatrix}
=
-\sum\limits_{r=1}^{m_i}\begin{bmatrix}[\widehat N_{i,r}^n]^\top\boldsymbol\lambda_{i,r}^n 
+ 
[\widehat N_{i,r}^t]^\top\boldsymbol\lambda_{i,r}^t\\ 0\end{bmatrix} =
-\begin{bmatrix}\widehat N_i^\top \boldsymbol\lambda_i + 
\widehat N_{i,\rm bc}^\top \boldsymbol\lambda_{i,\rm bc}
\\ 0\end{bmatrix},
%\qquad
%K_i\in\mathbb R^{(n_{i,\rm u}+n_{i,\rm p})\times(n_{i,\rm u}+n_{i,\rm p})},
\end{equation}
which is the counterpart of \eqref{eq:stokes_local_matrix} in the DD-PNM formulation. Here, the coupling matrices $\widehat N_i$ and $\widehat N_{i,\rm bc}$ are obtained by vertically stacking the per-interface blocks $\widehat N_{i,r}^n$ and $\widehat N_{i,r}^t$, and consistently for the vectors $\boldsymbol\lambda_{i}, \boldsymbol\lambda_{i,\rm bc}$. Moreover, we explicitly separate the contributions from prescribed inlet/outlet pressure tractions in $\boldsymbol\lambda_{i,\mathrm{bc}}$ and the associated coupling matrix $\widehat N_{i,\mathrm{bc}}$. 
% We stack Neumann coupling blocks vertically, first interface-wise, then normal and tangential in $\widehat N_i^\top$, 
Let $H_i:=[K_i^{-1}]_{11}\in\mathbb R^{n_{i,\rm u}\times n_{i,\rm u}}$ be the velocity--velocity block of $K_i^{-1}$. The subdomain velocity solution is thus
$$
\mathbf U_i = -H_i\widehat N_i^\top \boldsymbol\lambda_i - H_i\widehat N_{i,\rm bc}^\top \boldsymbol\lambda_{i,\rm bc}. 
$$
Velocity continuity across each physical interface $\omega$ implies
$$
\sum_{i:\omega\subset\partial\Omega_{i}}
\mathbf{U}_i|_\omega \cdot \mathbf{n}_{i} = 0,
\quad
\sum_{i:\omega\subset\partial\Omega_{i}}
\mathbf{U}_i|_\omega \cdot \mathbf{t}_{i} = 0,
$$
%which is 
%$$\sum_{i:\omega\subset\partial\Omega_{i}}(\widehat{N}_{i,\omega}^{n}\mathbf{U}_{i})=0,\quad\sum_{i:\omega\subset\partial\Omega_{i}}(\widehat{N}_{i,\omega}^{t}\mathbf{U}_{i})=0$$
Expressed in weak form, this yields the interface conditions
$$\sum_{i=1}^N(\widehat{R}_i^{\rm u})^\top\widehat{N}_i\:\mathbf{U}_i\:=\:0,$$
where $\widehat{R}_i^{\rm u}$ maps the stacked local internal-face DOFs of $\Omega_i$ to the global internal-face dof vector. Similarly, we introduce  $\widehat{R}_i^{\rm k}$ that maps local boundary-face DOFs of $\Omega_i$ 
 into the global boundary traction vector $\boldsymbol{\lambda}_{\rm bc}$.

Substituting the expression for $\mathbf{U}_i$ and grouping terms, we obtain the assembled C-DD system
$$\underbrace{\left(\sum_i(\widehat{R}_i^{\rm u})^\top \widehat{N}_i H_i \widehat{N}_i^\top \widehat{R}_i^{\rm u}\right)}_{\widehat{S}}\boldsymbol{\lambda}=\underbrace{-\left(\sum_i(\widehat{R}_i^{\rm u})^\top \widehat{N}_i H_i \widehat{N}_{i,\rm bc}^\top\widehat{R}_i^{\rm k}\right)\widehat{\boldsymbol\lambda}_\mathrm{bc}}_{\widehat{\mathbf F}},$$
which is the counterpart of the Schur complement system in DD-PNM formulation \eqref{eq:global_schur}. 

%with the local (fine) Schur/Steklov map
Here we state, without proof, the C-DD counterparts of the DD-PNM results: 
(i) the C-DD Schur operator  $\hat{S}$ is SPD; (ii) the energy identity 
\begin{equation}
\sum_{i=1}^N\|\mathbf{U}_i({\boldsymbol\lambda})-\mathbf{U}_i({\boldsymbol\mu})\|_{A_i}^2\:=\:\|{\boldsymbol\lambda}-{\boldsymbol\mu}\|_{\widehat{S}}^2 
\label{eq:energy_identity_c-DD}
\end{equation}
holds for any two internal tractions ${\boldsymbol\lambda},{\boldsymbol\mu}$, where 
the Schur energy norm and the subdomain viscous energies are defined as
$$\|\mathbf z\|_{\widehat{S}}^2:=\mathbf z^\top\widehat{S}\mathbf z,\quad\|\mathbf{U}_i\|_{A_i}^2:=\mathbf{U}_i^{\top}A_i\mathbf{U}_i.$$

%\begin{remark}
%
%
%\end{remark}

\subsection{DD-PNM formulation}
We now rederive the DD-PNM formulation as a restriction of the C-DD formulation introduced in the previous subsection.  
%let $r_i:=|\mathcal F_i^{\mathrm{int}}|$ be the number of local internal faces. 
Let $\boldsymbol\alpha\in\mathbb{R}^m$ collect one scalar per global internal face. For each local interface $\Gamma_{i,r}$ of $\Omega_i$, we define its normal projection $P_{i,r}\in\mathbb{R}^{m_{i,r}\times1}$ as the $L^2$ projection of the constant function $1$ onto the face basis, 
i.e., $M_{i,r}P_{i,r}=\mathbf b_{i,r}$ with $(M_{i,r})_{kl}=\int_{\Gamma_{i,r}}\widehat\varphi_{i,r,k}\widehat\varphi_{i,r,l}\,dx$ and
$(\mathbf b_{i,r})_k=\int_{\Gamma_{i,r}}\widehat\varphi_{i,r,k}\,dx$.
Then the fine (stacked) internal traction vector of $\Omega_i$ induced by the global vector $\boldsymbol\alpha$ is
$$\boldsymbol{\lambda}_i(\alpha)\:=\:\begin{bmatrix}\boldsymbol{\lambda}_i^n(\boldsymbol\alpha)\\\boldsymbol{\lambda}_i^t(\boldsymbol\alpha)\end{bmatrix}\:=\:\begin{bmatrix}P_i\,E_i\,\boldsymbol\alpha\\0\end{bmatrix},\quad P_i:=\mathrm{diag}\{P_{i,r}\}.$$
where $E_i$ is restriction matrix
mapping global interface indices to the local interface indices of $\Omega_i$. 
Assembling these local restrictions using $\widehat{R}_i^{\rm u}$, the global internal interface traction is
$$\boldsymbol{\lambda}=P\boldsymbol\alpha,\quad P:=\sum_i\widehat{R}_i^{\rm u}\begin{bmatrix}P_iE_i\\0\end{bmatrix}.$$
Insert $\boldsymbol{\lambda}_{\rm u}=P\boldsymbol\alpha$ into the classical DD system yields the DD-PNM system as a Galerkin restriction:
$$
\quad S\operatorname{\boldsymbol\alpha}=\mathbf g,\quad S:=P^{\top}\widehat{S}P,\quad \mathbf g:=P^{\top}\widehat{\mathbf F}.
$$

\begin{theorem}[Error estimate for DD--PNM]\label{thm:ddpnm-error}
Let $\boldsymbol\lambda^{\emph{FE}}$ be the classical DD  interface traction, i.e.
$\widehat S\,\boldsymbol\lambda^{\emph{FE}}=\widehat{\mathbf F}$.
Let $\boldsymbol\alpha^{\emph{DD-PNM}}$ be the solution of DD-PNM, i.e., the solution of
$S\,\boldsymbol\alpha=\mathbf g$, and set $\boldsymbol\lambda^{\emph{DD-PNM}}:=P\,\boldsymbol\alpha^{\emph{DD-PNM}}$.
Then:
\begin{equation}\label{eq:proj-identity}
\|\boldsymbol\lambda^{\emph{FE}}-\boldsymbol\lambda^{\emph{DD-PNM}}\|_{\widehat S}
\ =\
\min_{\boldsymbol\alpha\in\mathbb R^{m}} \ \|\boldsymbol\lambda^{\emph{FE}}-P\boldsymbol\alpha\|_{\widehat S},
\end{equation}
and the energy identity
\begin{equation}\label{eq:field-equality_error_estimate}
\sum_{i=1}^N \|\mathbf U_i^{\emph{FE}}-\mathbf U_i^{\emph{DD-PNM}}\|_{A_i}^{2}
\ =\
\|\boldsymbol\lambda^{\emph{FE}}-\boldsymbol\lambda^{\emph{DD-PNM}}\|_{\widehat S}^{2}
\ =\
\min_{\boldsymbol\alpha}\ \|\boldsymbol\lambda^{\emph{FE}}-P\boldsymbol\alpha\|_{\widehat S}^{2}.
\end{equation}
\end{theorem}

\begin{proof}
\emph{(i) } 
Subtracting the C-DD and DD-PNM systems
 gives the Galerkin orthogonality
\[
P^\top \widehat S\big(\boldsymbol\lambda^{\mathrm{DD-PNM}}-\boldsymbol\lambda^{\mathrm{FE}}\big)=0,
\]
Thus the error $e:=\boldsymbol\lambda^{\mathrm{FE}}-\boldsymbol\lambda^{\mathrm{DD-PNM}}$ is orthogonal to the
restriction space $\mathcal R(P)$ with respect to the inner product
$\langle \mathbf x,\mathbf y\rangle_{\widehat S}:=\mathbf x^\top \widehat S\mathbf y$. For any $\boldsymbol\alpha$, applying Pythagoras identity gives
\[
\|\boldsymbol\lambda^{\mathrm{FE}}-P\boldsymbol\alpha\|_{\widehat S}^2
= \|e\|_{\widehat S}^2 + \|\boldsymbol\lambda^{\mathrm{DD-PNM}}-P\boldsymbol\alpha\|_{\widehat S}^2
\ \ge\ \|e\|_{\widehat S}^2,
\]
with equality if and only if $\boldsymbol\alpha=\boldsymbol\alpha^{\mathrm{DD-PNM}}$. This proves \eqref{eq:proj-identity}.

\emph{(ii)} 
the second equality follows by apply the energy identity \eqref{eq:energy_identity_c-DD} with $\boldsymbol\lambda=\boldsymbol\lambda^{\mathrm{FE}}$ and $\boldsymbol\mu=\boldsymbol\lambda^{\mathrm{DD-PNM}}$.
\end{proof}

\begin{remark}
The previous theorem shows that the DD–PNM approximation error equals the minimal distance, in Schur complement energy, between the FE traction and the DD–PNM space $\mathcal R(P)$. This error typically does not diminish under mesh refinement unless the interface approximation space is suitably enriched (for example, by adding higher-order face moments or tangential modes). 
\label{remark:error_estimate}
\end{remark}

%Consequently, if the FE normal traction is not facewise constant on some interfaces in the $h\!\downarrow\!0$ limit,
%the modeling distance $\mathrm{dist}_S(\bar\alpha,\cdot)$ does not vanish and a refinement plateau occurs unless the interface model space is enriched (e.g., add face moments).

%\emph{DD--PNM is the $\widehat S_{uu}$-orthogonal projection of the full FE traction onto the subspace
%$\{[P_nDp;\ 0]\}$.}

%\bibliographystyle{plain}	

\bibliographystyle{siamplain}
\bibliography{all_references}
\end{document}